\documentclass[10pt]{article}

\usepackage{amsmath}
\usepackage{amsthm}
\usepackage{amssymb}
\usepackage{graphicx}
\usepackage{hyperref}

\newtheorem{thm}{Theorem}

\newtheorem{remark}{Remark}

\newtheorem{lemma}[thm]{Lemma}
\newtheorem{prop}[thm]{Proposition}
\newtheorem{cor}[thm]{Corollary}

\newcommand{\cO}{c_0}

\newcommand{\grad}{{\rm grad\,}}

\newcommand{\I}[1]{\left\|#1\right\|}

\newcommand{\ip}[2]{\left\langle #1, #2 \right\rangle}

\newcommand{\R}{{\mathbb R}}
\newcommand{\Rp}{{\mathbb R}_+}
\newcommand{\lb}{\left(}
\newcommand{\rb}{\right)}
\newcommand{\lsb}{\left[}
\newcommand{\rsb}{\right]}
\newcommand{\p}{\partial}

\renewcommand{\O}[1]{\mathrm{O}\lb #1\rb}

\newcommand{\lra}[1]{\langle#1\rangle}

\renewcommand{\L}{{\mathcal L}}

\title{Blow-up in Nonlinear Heat Equations}
\begin{document}
\author{S. Dejak$^{1}$$^*$, Zhou Gang$^{1}$\thanks{Supported by NSERC under Grant NA7901.}, I.M.Sigal$^{1}$$^*$, S. Wang$^{2}$}
\maketitle 
\centerline{\small{$^{1}$Department of Mathematics, University of
Toronto, Toronto, Canada}} \centerline{\small{$^{2}$Department of
Mathematics, University of Notre Dame, Notre Dame, U.S.A.}}
\section*{Abstract}
In this paper we study the blowup problem of nonlinear heat equations. Our result show that for a certain family of initial conditions the solution will blowup in finite time, the blowup parameters satisfy some dynamics which are asymptotic stable, moreover we provide the remainder estimates. Compare to the previous works our approach is analogous to one used in bifurcation
theory and our techniques can be regarded as a time-dependent
version of the Lyapunov-Schmidt decomposition.
\section{Introduction}
\label{Intro} We study the blow-up problem for the one-dimensional
nonlinear heat equations (or the reaction-diffusion equations) of
the form
\begin{equation}\label{NLH}
\begin{array}{lll}
u_{t}&=&\partial_{x}^{2} u+|u|^{p-1}u\\
u(x,0)&=&u_{0}(x)
\end{array}
\end{equation} with $p>1$.
Equation \eqref{NLH} arises in the problem of heat flow and the
theory of chemical reactions. Similar equations appear in the motion
by mean curvature flow (see \cite{SS}), vortex dynamics in
superconductors (see \cite{CHO,MZ2}), surface diffusion (see
\cite{BBW}) and chemotaxis (see \cite{BCKSV,BB}). Equation
\eqref{NLH} has the following properties:
\begin{itemize}
 \item \eqref{NLH} is invariant
 with respect to the scaling transformation,
 \begin{equation}\label{rescale}
 u(x,t)\rightarrow \lambda^{\frac{2}{p-1}}
u(\lambda x,\lambda^{2} t)
\end{equation} for any constant $\lambda>0,$ i.e. if
$u(x,t)$ is a solution, so is $\lambda^{\frac{2}{p-1}}u(\lambda
x,\lambda^{2}t).$
 \item \eqref{NLH} has $x-$independent of $x$ (homogeneous)
 solutions:
\begin{equation}
 u_{hom}=[u_{0}^{-p+1}-(p-1)t]^{-\frac{1}{p-1}}.
 \label{eqn:3}
\end{equation} These solutions blow up in finite time
 $t^* = \lb (p-1)u_0^{p-1}\rb^{-1}$ for $p>1$.
\item \eqref{NLH} is an $L^2$-gradient system $\p_t u = -\grad {\cal E}(u)$,
with the energy
\begin{equation}
{\cal E}(u):= \int \frac{1}{2} u_x^2-\frac{1}{p+1}u^{p+1}.
\label{eqn:L2Energy}
\end{equation}
\end{itemize}
(With the $L^2(\R)$ metric, $\grad {\cal E}$ is defined by the
relation $\p{\cal E}(u)\xi=\ip{\grad{\cal E}(u)}{\xi}$, so that $
\grad {\cal E}(u)=-(\p_x^2 u+u^p).) $ We immediately have that the
energy ${\cal E}$ decreases under the flow of \eqref{NLH}.

The linearization of \eqref{NLH} around $u_{hom}$ shows that the
solution $u_{hom}$ is unstable.  Moreover, it is shown in \cite{GK1}
that if either $n\le 2$ or $p\le (n+2)/(n-2)$, then \eqref{NLH} in
dimension $n$ has no other self-similar solutions of the form
$(T-t)^{-\frac{1}{p-1}}\phi\lb x/\sqrt{T-t}\rb$, $\phi\in
L^{\infty}$, besides $u_{hom}$.

The local well-posedness of \eqref{NLH} is well known (see, e.g.
\cite{Ball} for $H^\alpha$, $0\le\alpha<2$). Moreover for some data
$u_{0}(x)$, the solutions $u(x,t)$ might blowup in finite time
$T>0$. Thus, two key problems about \eqref{NLH} are
\begin{enumerate}
 \item Describe initial conditions for which solutions of Equation
\eqref{NLH} blowup in finite time;
 \item Describe the blowup profile of such solutions.
\end{enumerate}

It is expected (see e.g. \cite{BrKu}) that the blowup profile is
universal $-$ it is independent of lower power perturbations of the
nonlinearity and of initial conditions within certain spaces.

There is rich literature regarding the blowup problem for Equation
(~\ref{NLH}). We review quickly relevant results.  Starting with
\cite{Fu}, various criteria for blow-up in finite time were derived,
see e.g. \cite{Fu,Ball,CHI,EVA,Le1,Le2,OHLM,Qu,Sou,MR1843848, MR1433084}. For example, if
$u_{0}\in \mathcal{H}^{1}\cap L^{p+1}$ and $\mathcal{E}(u_{0})<0$,
where ${\cal E}(u)$ is the energy functional for \eqref{NLH} defined
in \eqref{eqn:L2Energy}, then it is proved in \cite{Le1} that
$\|u(t)\|_2^{2}$ blows up in finite time $t^*$. By the observation
$$\frac{1}{2}\frac{d}{dt}\|u(t)\|_2^{2}\leq
\|u(t)\|_\infty^{p-1}\|u(t)\|_2^{2}$$ we have that $\|u(t)\|_\infty$
blows up in finite time $t^{**}\leq t^*$ also.  (In this paper, we
denote the norms in the $L^p$ spaces by $\|\cdot\|_p$.)

Recall that a solution $u(x,t)$ is said to blowup at time $t^*$ if
it exists in $L^\infty$ for $[0,t^*)$ and
$\sup_{x}|u(x,t)|\rightarrow\infty$ as $t\rightarrow t^*$.  The
first result on asymptotics of the blowup was obtained in the
pioneering paper \cite{GK1} where the authors show that under the
conditions
\begin{equation}\label{eq:boundedness}
|u(x,t)|(t_*-t)^\frac{1}{p-1}\ \mbox{is bounded on}\ B_1\times
(0,t_*),
\end{equation}
where $B_1$ is the unit ball in $\R^n$ centred at the origin, and
either $p\le\frac{n+2}{n-2}$ or $n\le 2$ and assuming blowup takes
place at $x=0$, one has
\begin{equation*}
\displaystyle\lim_{\lambda\rightarrow 0} \lambda^\frac{2}{p-1}
u(\lambda x,t_{*}+\lambda^2(t-t_*) )=\pm
\lb\frac{1}{p-1}\rb^\frac{1}{p-1}(t_*-t)^{-\frac{1}{p-1}}\
\mbox{or}\ 0.
\end{equation*}
This result was further improved in several papers (see e.g.
\cite{GK2, GK3, HV92, MR1164066, Mer,Vel,MR1230711,
MR1317705,MR1433084,BrKu, MZ3, MZ4, MZ5}). A blowup solution
satisfying the bound (~\ref{eq:boundedness}) is said to be of type
I. This bound was proven under various conditions in \cite{GK2, MZ3,
MZ4, Wei2, GiMaSa}.  Furthermore, the limits of $H^1$-blowup
solutions $u(x,t)$ as $t\uparrow T$, outside the blowup sets were
established in \cite{HV92,MR1164066, Mer, Vel,MR1230711,
MR1317705,MR1433084,BrKu, MZ5, FMZ}.

For $p>1$, Herrero and Vel\'{a}zquez \cite{HV2} (see also
~\cite{MR1230711}) proved that if the initial condition $u_0$ is
continuous, nonnegative, bounded, even and has only one local
maximum at $0$, and if the corresponding solution blows up, then
\begin{equation}\label{eq:asy}
\lim_{t\uparrow t^*}(t^*-t)^{\frac{1}{p-1}}u(y((t^*-t)ln|
t^*-t|)^{1/2},t)=(p-1)^{-\frac{1}{p-1}}[1+\frac{p-1}{4p}y^{2}]^{-\frac{1}{p-1}}
\end{equation} uniformly on sets $|y|\leq R$ with $R>0$.  Further
extensions of this result are achieved in \cite{HV92, Vel, MR1164066, MR1230711}.

Later Bricmont and Kupiainen \cite{BrKu} constructed a co-dimension
2 submanifold, of initial conditions such that \eqref{eq:asy} is
satisfied on the whole domain. More specifically, given a small
function $g$ and a small constant $b>0$, they find constants $d_{0}$
and $d_{1}$ depending on $g$ and $b$ such that the solution to
\eqref{NLH} with the datum
\begin{equation}\label{dis}
u_{0}^{*}(x)
=(p-1+bx^{2})^{-\frac{1}{p-1}}(1+\frac{d_{0}+d_{1}x}{p-1+bx^{2}})^{\frac{1}{p-1}}+g(x)\end{equation}
has the convergence \eqref{eq:asy} uniformly in
$y\in(-\infty,+\infty)$. The result of \cite{BrKu} was generalized
in \cite{MZ1,FMZ} (see also \cite{GaPo1986}), where it is shown that
there exists a neighborhood $U$, in the space
$\mathcal{H}:=L^{p+1}\cap \mathcal{H}^{1}$, of $u_{0}^{*}$, given in
\eqref{dis}, such that if $u_{0}\in U$, then the solution $u(x,t)$
blows up in a finite time $t^*$ and satisfies (~\ref{eq:asy}) for
$x\in \mathbb{R}$. They conjectured that this asymptotic behavior is
generic for any blow-up solution.

The starting point in the above works, which goes back to Giga and
Kohn \cite{GK1}, is passing to the similarity variables
$y:=x/\sqrt{t^*-t}$ and $s:=-\log (t^*-t)$, where $t^*$ is the
blowup time, and to the rescaled function
$w(y,s)=(t^*-t)^\frac{1}{p-1} u(x,t)$.  Then one studies the
resulting equation for $w$:
\begin{equation}\label{eq:freDir}
\p_s w=\p_y^2 w-\frac{1}{2} y\p_y w-\frac{1}{p-1} w+|w|^{p-1}w.
\end{equation}
Most of the work above uses relations involving the energy
functional
\begin{equation}\label{eq:energy}
S(w):=\frac{1}{2}   \int \lb |\nabla w|^2 +\frac{1}{p-1} |w|^2
-\frac{2}{p+1} |w|^{p+1}\rb e^{-\frac{1}{4} y^2}\, dy,
\end{equation}
introduced in \cite{GK1}, and related functionals. In particular,
one uses the relation
\begin{equation}
\p_s S(w)=-  \int |\p_s w|^2 e^{-\frac{1}{4} y^2}\, dy. \label{eqn:9a}
\end{equation}
\begin{remark}
Equation (~\ref{eq:freDir}) is the gradient system $\p_s w=-\grad
S(u)$ in the metric space ${L^2(e^{-\frac{a}{4} y^2}\,dx)}$. ($\grad
S(u)$ is defined by the equation $ \p S(u)\xi=\ip{\grad
S(u)}{\xi}_{L^2(e^{-\frac{a}{4} y^2}\, dy)}$.) Hence $S$ decreases
under the flow of \eqref{NLH} and so \eqref{eqn:9a} implies that
$\p_s w\rightarrow 0$ as $s\rightarrow \infty$.
\end{remark}
Blowup as a single point was studied as early as \cite{Wei1} (see
also ~\cite{MR1433084}). In 1992, Merle \cite{Mer} proved that given
an finite number of points $x_1$, $x_2$, $\ldots$, $x_k$ in
$I=(-1,1)$ (or any other domain $I$ in $\R$), there is a positive
solution to the nonlinear heat equation which blowups up at time $T$
with blowup points $x_1$, $x_2$, $\ldots$, $x_k$.  This theorem can
be generalized to allow the sign ($+\infty$ or $-\infty$) to be
chosen at each blowup point $x_i$.

In this paper, we consider \eqref{NLH} with initial conditions which
are even, have, modulo a small perturbation, a maximum at the
origin, are slowly varying near the origin and are sufficiently
small, but not necessarily vanishing, for large $|x|$. In
particular, the energy ${\cal E}(u)$ for such initial conditions
might be infinite. We show that the solutions of \eqref{NLH} for
such initial conditions blowup in a finite time $t^*$ and we
characterize asymptotic dynamics of these solutions. As it turns
out, the leading term is given by the expression
\begin{equation}
\lambda(t)^\frac{2}{p-1}\lsb\frac{2 c(t)}{p-1+b(t)\lambda(t)^2
x^2}\rsb^\frac{1}{p-1} \label{leading}
\end{equation}
(cf \eqref{eq:asy}) where the parameters $\lambda(t)$, $b(t)$ and
$c(t)$ obey certain dynamical equations whose solutions give
\begin{equation}\label{eq:blowdy}
\begin{array}{lll}
\lambda(t)&=\lambda_0(t^*-t)^{-\frac{1}{2}}(1+o(1))\\
b(t)&=\frac{(p-1)^{2}}{4p|ln|t^*-t||}(1+O(\frac{1}{|ln|t^*-t||^{1/2}}))\\
c(t)&=\frac{1}{2}-\frac{p-1}{4p|ln|t^*-t||}(1+O(\frac{1}{ln|t^*-t|})).
\end{array}
\end{equation}
with $\lambda(0)=\sqrt{2 \cO+\frac{2}{p-1} b_0}$, $c_0,\ b_{0}>0$
depends on the initial datum. Here $o(1)$ is in $t^*-t$. Moreover,
we estimate the remainder, the difference between $u(x,t)$ and
\eqref{leading}. Our techniques are different from the papers
mentioned above, the closest to our approach is \cite{BrKu}. Our
main point is that we do not fix the time-dependent scale in the
self-similarity (blowup) variables but let its behaviour, as well as
behaviour of other parameters ($b$ and $c$) to be determined by the
equation.  This approach is analogous to one used in bifurcation
theory and our techniques can be regarded as a time-dependent
version of the Lyapunov-Schmidt decomposition.

In what follows we use the notation $f\lesssim g$ for two functions
$f$ and $g$ satisfying $f\leq Cg$ for some universal constant $C$.
We will also deal, without specifying it, with weak solutions of
Equation (~\ref{NLH}) in some appropriate sense. These solutions can
be shown to be classical for $t>0.$
\begin{thm}\label{maintheorem} Suppose in \eqref{NLH} the initial datum $u_{0}\in L^\infty(\R)$ is even
and satisfy
\begin{equation}\label{INI2}
\|\langle
x\rangle^{-n}(u_{0}(x)-(\frac{2\cO}{p-1+b_{0}x^{2}})^{\frac{1}{p-1}})\|_{\infty}\le
\delta_{n}
\end{equation} with $n= 0, 3, 1/2\leq \cO\leq 2$, $0\leq b_{0}, \delta_{0}\ll 1$
and $\delta_3 = C b_{0}^2$. Then
\begin{enumerate}
\item[(1)] There exists a time $t^*\in (0,\infty)$ such that the
solution $u(x,t)$ blows up at $t\rightarrow t^{*}.$
\item[(2)] When
$t\leq t^*$, there exist unique positive, $C^{1}$ functions
$\lambda(t)$, $b(t)$ and $c(t)$ with $b(t)\lesssim b_0$ such that
$u(x,t)$ can be decomposed as
$$u(x,t)=\lambda^{\frac{2}{p-1}}(t)[(\frac{2c(t)}{p-1+b(t)\lambda^{2}(t)x^{2}})^{\frac{1}{p-1}}+\eta(x,t)]$$
with the fluctuation part, $\eta,$ admitting the estimate $\|\langle
\lambda(t)x\rangle^{-3}\eta(x,t)\|_{\infty}\lesssim b^{2}(t).$
\item[(3)] The functions $\lambda(t)$, $a(t)$, $b(t)$ and $c(t)$
are of the form (~\ref{eq:blowdy})..
\end{enumerate}
\end{thm}
The proof is given in Section ~\ref{SecMain}. Thus our result shows
the blow-up at $0$ for a certain neighborhood of the homogeneous
solution, \eqref{eqn:3}, with a detailed description of the leading
term and an estimate of the remainder in $L^\infty$. In fact, we
have not only the asymptotic expressions for the parameters $b$ and
$c$ determining the leading term and the size of the remainder, but
also dynamical equations for these parameters:
\begin{align}
 b_\tau&=-\frac{4 p}{(p-1)^{2}}
b^2+c^{-1}c_{\tau}
b+{\mathcal R}_b(\eta,b,c),\\
c^{-1}c_\tau&=2(\frac{1}{2}-c)-\frac{2}{p-1}b+{\mathcal
R}_c(\eta,b,c),
\end{align}
where $\tau$ is a 'blow-up' time related to the original time $t$ as
$\tau(t):=  \int_{0}^{t}\lambda^{2}(s)ds$, the remainders have the
estimates
\begin{equation}
\begin{array}{lll}
& &{\mathcal R}_b(\eta,b,c),{\mathcal R}_c(\eta,b,c)\\
&=&\O{ b^3+[|c-\frac{1}{2}|+|c_\tau|]b^2+|b_\tau|
b+b\|\eta(\cdot,t)\|_{X}+\|\eta(\cdot,t)\|_{X}^2+\|\eta(\cdot,t)\|_{X}^p}.
\end{array}
\end{equation}
with the norm $\|\eta(\cdot,t)\|_{X}:=\|\langle
\lambda(t)x\rangle^{-3}\eta\|_{\infty}.$
\begin{remark}
\begin{enumerate}
\item[(a)]
The restriction \eqref{INI2} on the initial condition $u_0(x)$
states roughly that $mod\ O(b^{2})$ $u_0(x)$ (after initial
rescaling if necessary) has a form $\phi(\sqrt{b(0)} x)$ for
$|x|\lesssim \frac{1}{\sqrt{b}}$ with an absolute maximum at $x=0$
and is of the size $\delta_0$ for $|x|\gg \frac{1}{\sqrt{b}}$.
\item[(b)] We allow for initial conditions to have infinite energy. It seems that previously, blowup for the nonlinear heat
equation was studied only for finite energy solutions.
\item[(c)] We expect our approach can be extended to general data,
to more general nonlinearities and to dimensions $\geq 2.$
\end{enumerate}
\end{remark}
This paper is organized as follows.  In Sections
\ref{Section:BV}-\ref{Section:Reparam} we present some preliminary
derivations and some motivations for our analysis. In Section
\ref{Section:APriori}, we formulate a priori bounds on solutions to
\eqref{NLH} which are proven in Sections \ref{Section:9},
\ref{SEC:EstM1} and \ref{SEC:EstM2}. We use these bounds in Section
\ref{SecMain} in order to prove our main result, Theorem
\ref{maintheorem}. In Sections \ref{Section:Splitting},
\ref{Section:Rescaling} and \ref{Section:PropEst} we lay the ground
work for the proof of the a priori bounds of Section
\ref{Section:APriori}. In particular, in Section
\ref{Section:Splitting}, using a Lyapunov-Schmidt-type argument we
derive equations for the parameters $a$, $b$ and $c$ and fluctuation
$\eta$.  In Section \ref{Section:Rescaling} we rescale our equations
in a convenient way and in Section \ref{Section:PropEst} we estimate
the corresponding propagators.  As was mentioned above, the results
of Sections \ref{Section:Splitting}, \ref{Section:Rescaling} and
\ref{Section:PropEst} are used in Sections \ref{Section:9},
\ref{SEC:EstM1} and \ref{SEC:EstM2} in order to prove the a priori
estimates. The paper has four appendices. In Appendix
\ref{Appendix:LWP}, we present a local existence result for
\eqref{NLH} in the $L^\infty$ space and a blowup criterion.  In
Appendix \ref{sec:BlowUpDyn}, we discuss other relations between the
parameters $a$, $b$ and $c$ than the one used in the paper
($c=\frac{1}{2}a+\frac{1}{4}$).  In Appendix
\ref{Appendix:SpectrumLinear} we investigate the spectrum of the
linearized operator.  The result of this appendix is not used in the
main part of this paper. In Appendix \ref{Sec:Trotter}, we prove a
convenient form of the Feynmann-Kac-type formula.  It seems the
results of Appendices \ref{Appendix:LWP} and \ref{Sec:Trotter} are
generally assumed to be known, but we did not find them in the
literature, so we included them for the reader's convenience.
\section*{Acknowledgements}
We would like to thank V.S.Buslaev, C.Fefferman, G.M.Graf and
A.Soffer for helpful discussions and, especially, Kong Wenbin, Zou
Xiangqun for careful reading of the manuscript and many useful
suggestions.

\section{Blow-Up Variables and Almost Solutions}
\label{Section:BV}
In this section we pass from the original variables $x$ and $t$ to
the blowup variables $y:=\lambda(t) (x-x_{0}(t))$ and
$\tau:=  \int_{0}^{t} \lambda^2(s) \,ds$.  The point here is that we
do not fix $\lambda(t)$ and $x_0(t)$ but consider them as free
parameters to be found from the evolution of \eqref{NLH}. Assume for
simplicity that $0$ is a maximum point of $u_0$ and that $u_0$ is
even with respect to $x=0$. In this case $x_{0}$ can be taken to be
$0$. Suppose $u(x,t)$ is a solution to \eqref{NLH} with an initial
condition $u_{0}(x)$. We define the new function
\begin{equation}\label{eq:definev}
v(y,\tau):=\lambda^{-\frac{2}{p-1}}(t)u(x,t)
\end{equation} with $y:=\lambda(t)x$ and $\tau:= \int_{0}^{t}\lambda^{2}(s)ds.$
The function $v$ satisfies the equation
\begin{equation}\label{eqn:BVNLH}
v_{\tau}=\lb\p_y^2-a y\p_y-\frac{2 a}{p-1}\rb v+|v|^{p-1}v.
\end{equation}
where $a:=\lambda^{-3}\p_t\lambda$.  The initial condition is
$v(y,0)=\lambda^{-\frac{2}{p-1}}_{0} u_{0}( y/\lambda_{0})$, where $\lambda_0$ is an
initial condition for the scaling parameter $\lambda$.

If the parameter $a$ is a constant, then \eqref{eqn:BVNLH} has the
following homogeneous, static (i.e. $y$ and $\tau$-independent)
solutions
\begin{equation}
v_a:=\lb \frac{2 a}{p-1}\rb^\frac{1}{p-1}.
\end{equation}
In the original variables $t$ and $x$, this family of solutions
corresponds to the homogeneous solution \eqref{eqn:3} of the
nonlinear heat equation with the parabolic scaling $\lambda^{-2}=2
a(T-t)$, where the blowup time,
$T:=\left[u_0^{p-1}(p-1)\right]^{-1}$, is dependent on $u_0$, the
initial value of the homogeneous solution $u_{hom}(t)$.

If the parameter $a$ is $\tau$ dependent but $|a_\tau|$ is small,
then the above solutions are good approximations to the exact
solutions. Another approximation is the solution of the equation $ a
y v_y+\frac{2 a}{p-1} v=v^p$, obtained from \eqref{eqn:BVNLH} by
neglecting the $\tau$ derivative and second order derivative in $y$.
This equation has the general solution
\begin{equation}
v_{a b}:=\lb\frac{2 a}{p-1 + b y^2}\rb^\frac{1}{p-1}
\end{equation}
for all $b\in\R$.  (The above equation is equivalent to the equation
$ \p_y\lb y^\frac{2}{p-1}v \rb=\frac{1}{a y^3}\lb y^\frac{2}{p-1}
v\rb^p$). In what follows we take $b\ge 0$ so that $v_{a b}$ is
nonsingular.  Note that $v_{0a}=v_a$.


\section{"Gauge" Transform}
\label{sec:Gauge}
We assume that the parameter $a$ depends slowly on $\tau$ and treat
$|a_\tau|$ as a small parameter in a perturbation theory for
Equation \eqref{eqn:BVNLH}.
 In order to convert the global non-self-adjoint operator $a y\p_y$
appearing in this equation into a more tractable local and
self-adjoint operator we perform a gauge transform. Let
\begin{equation}\label{eq:definew}
w(y,\tau):=e^{-\frac{a y^2}{4}}v(y,\tau).
\end{equation}
Then $w$ satisfies the equation
\begin{equation} \label{eqn:w}
w_\tau=\lb \p_y^2-\frac{1}{4}\omega^2
y^2-\lb\frac{2}{p-1}-\frac{1}{2}\rb a \rb
w+e^{\frac{a}{4}(p-1)y^2}|w|^{p-1}w,
\end{equation}
where $\omega^2=a^2+a_\tau$.  The approximate solution $v_{a b}$ to
\eqref{eqn:BVNLH} transforms to $v_{a b c}$ where $v_{a b c}:=v_{c b
}e^{-\frac{a y^2}{4}}$. Explicitly
\begin{equation}
v_{a b c}:=\lb\frac{2 c}{p-1+b y^2}\rb^{\frac{1}{p-1}}e^{-\frac{a
y^2}{4}}.
\end{equation}

Equation \eqref{eqn:w} is the $L^2$-gradient system with the energy
\begin{equation}\label{eqn:energy}
{\cal E}(w):=- \int \frac{1}{2} w \lb \p_y^2-\frac{1}{4}\omega^2
y^2-\lb\frac{2}{p-1}-\frac{1}{2}\rb a \rb w-\frac{1}{p+1}
e^{\frac{a}{4}(p-1) y^2}|w|^{p+1}\, dy.
\end{equation}
This energy is related to the functional (~\ref{eq:energy}). It
satisfies the relation
\begin{equation*}
\p_\tau{\cal E }(w)(\tau)=-  \int |\p_s w|^2 e^{-\frac{1}{4} y^2}\,
dy.
\end{equation*}
Indeed, multiplying \eqref{eqn:w} by $w_\tau$, integrating over
space and then using that the linear operator in (~\ref{eqn:w}) is
self-adjoint gives this relation.
\section{Reparametrization of Solutions}
\label{Section:Reparam}
In this section we split solutions to \eqref{eqn:w} into the leading
term - the almost solution $v_{a b c}$ - and a fluctuation $\xi$
around it.  More precisely, we would like to parametrize a solution
by a point on the manifold $M_{as}:=\{v_{a b c}\, |\, a,b,c\in \Rp,
b \leq \epsilon,\, a=a(b,c)\}$ of almost solutions and the
fluctuation orthogonal to this manifold (large slow moving and small
fast moving parts of the solution).  Here $a=a(b,c)$ is a twice
differentiable function of $b$ and $c$.  For technical reasons, it
is more convenient to require the fluctuation to be almost
orthogonal to the manifold $M_{as}$. More precisely, we require
$\xi$ to be orthogonal to the vectors $\phi_{0 a}:=e^{-\frac{a}{4}
y^2}$ and $\phi_{2 a}:=(1-a y^2)e^{-\frac{a}{4} y^2}$ which are
almost tangent vectors to the above manifold, provided $b$ is
sufficiently small.  Note that $\xi$ is already orthogonal to
$\phi_{1 a}:=\sqrt{a} y e^{-\frac{a}{4} y^2}$ since our initial
conditions, and therefore, the solutions are even in $x$.

In this section and the rest of the paper except Appendix
~\ref{sec:BlowUpDyn} we fix the relation between the parameters $a$,
$b$ and $c$ as
$$2c=a+\frac{1}{2}.$$  In Appendix ~\ref{sec:BlowUpDyn}
we prove that under some conditions different functions of
$a=a(c,b)$ can be used.

Let $V_{a b}:=(\frac{2c}{p-1+by^{2}})^{\frac{1}{p-1}}$ with
$c=\frac{1}{2}a+\frac{1}{4}$.  We define a neighborhood:
\begin{equation*}
U_{\epsilon_0}:=\{v\in L^\infty(\R)\ |\ \|e^{-\frac{1}{9}
y^2}(v-V_{ab})\|_\infty=o(b)\ \mbox{for some}\ a\in[1/4,1],\
b\in[0,\epsilon_0]\ \}.
\end{equation*}

\begin{prop}\label{Prop:Splitting}
There exist an $\epsilon_{0}>0$ and a unique $C^1$ functional
$g:U_{\epsilon_0}\rightarrow \mathbb{R}^{+}\times \mathbb{R}^{+}$,
such that any function $v\in U_{\epsilon_0}$ can be uniquely written
in the form
\begin{equation} \label{eqn:split}
v =V_{g(v)} + \eta,
\end{equation}
with $\eta\perp e^{-\frac{a}{4} y^2}\phi_{0 a},\ e^{-\frac{a}{4}
y^2}\phi_{2 a}$ in $L^2(\R)$, $(a,b)=g(v)$. Moreover, if $(a_0,
b_0)\in [\frac{1}{4},1]\times [0,\varepsilon_0]$ and
$\|e^{-\frac{1}{9}y^{2}}(v-V_{a_{0},b_{0}})\|_{\infty}=o(b_{0})$,
then
\begin{equation}
\label{eqn:28a} |g(v)-(a_{0},b_{0})|\lesssim \|e^{-\frac{1}{9}a
y^2}(v-V_{a_{0}b_{0}})\|_\infty.
\end{equation}
\end{prop}
\begin{proof}
 The orthogonality conditions on the fluctuation can be written as
$G(\mu,v)=0$, where $\mu=(a,b)$ and $G:\mathbb{R}^{+}\times
\mathbb{R}^{+}\times
 L^\infty\lb \R\rb\rightarrow \R^2$ is defined
as
\begin{equation*}
G(\mu, v):=\lb \begin{array}{c} \ip{V_{\mu}-v}{e^{-\frac{ay^{2}}{4}}\phi_{0a}}\\
\ip{V_{\mu}-v}{e^{-\frac{ay^{2}}{4}}\phi_{2a}}
\end{array} \rb.
\end{equation*}
Here and in what follows, all inner products are $L^2$ inner
products.  Let $X:= e^{\frac{1}{9} y^2} L^\infty$ with the
corresponding norm.  Using the implicit function theorem we will
prove that for any $\mu_0:=(a_0, b_0)\in [\frac{1}{4},1]\times
(0,\epsilon_{0})$ there exists a unique $C^1$ function
$g:X\rightarrow \mathbb{R}^{+}\times \mathbb{R}^{+}$ defined in a
neighborhood $U_{\mu_0}\subset X$ of $V_{\mu_0}$ such that
$G(g(v),v)=0$ for all $v\in U_{\mu_0}$. Let
$B_{\varepsilon}(V_{\mu_0})$ and $B_\delta(\mu_0)$ be the balls in
$X$ and $\mathbb{R}^2$ around $V_{\mu_0}$ and $\mu_0$ and of the
radii $\varepsilon$ and $\delta$, respectively.

Note first that the mapping $G$ is $C^1$ and $G(\mu_0, V_{\mu_0})=0$
for all $\mu_0$.  We claim that the linear map $\p_\mu G(\mu_0,
V_{\mu_0})$ is invertible.  Indeed, we compute
\begin{equation}\label{linearization}
\p_\mu G(\mu,
v)|_{\mu=\mu_{0}}=A_{1}(\mu)+A_{2}(\mu,v)|_{\mu=\mu_{0}}
\end{equation}
where
$$
A_{1}(\mu):=\lb
\begin{array}{cc}
\ip{\p_a V_{
\mu}}{ e^{-\frac{a}{2} y^2}} & \ip{\p_b V_{\mu}}{ e^{-\frac{a}{2} y^2}}  \\
\ip{\p_a V_{\mu}}{(1-a y^2)e^{-\frac{a}{2} y^2}}& \ip{\p_b
V_{\mu}}{(1-a y^2)e^{-\frac{a}{2} y^2}}
\end{array}
\rb
$$ and
\begin{equation*}
A_{2}(\mu,v):= -\frac{1}{4}\lb\begin{array}{cc}
\langle V_{\mu}-v, y^{2}e^{-\frac{a}{2}y^{2}}\rangle& 0 \\
\ip{V_{\mu}-v}{\lb 1-a y^2\rb y^2 e^{-\frac{a}{2} y^2}} & 0
\end{array}\rb.
\end{equation*}
For $b>0$ and small, we expand the matrix $A_{1}$ in $b$ to get
$A_{1}=G_{1}G_{2}+O(b)$, where the matrices $G_{1}$ and $G_{2}$ are
defined as
$$G_{1}:=\left(
\begin{array}{ccc}
\langle -y^{2}e^{-\frac{ay^{2}}{4}},e^{-\frac{ay^{2}}{4}}\rangle &
\frac{1}{a+\frac{1}{2}}\langle e^{-\frac{ay^{2}}{4}},e^{-\frac{ay^{2}}{4}}\rangle\\
\langle -y^{2}e^{-\frac{ay^{2}}{4}},
(1-ay^{2})e^{-\frac{ay^{2}}{4}}\rangle & 0
\end{array}
\right)$$ and
$$G_{2}:=(\frac{a+1/2}{p-1})^{\frac{1}{p-1}}\frac{1}{p-1}\left(
\begin{array}{lll}
\frac{p-1}{4} & 1\\
1 & 0
\end{array}
\right).$$ Obviously the matrices $G_{1}$ and $G_{2}$ have uniformly
(in $a\in [\frac{1}{4},1]$) bounded inverses. Furthermore, by the
Schwarz inequality
$$\|A_{2}(\mu_{0},v)\|\lesssim \| v-V_{a_{0}b_{0}}\|_X.$$
Therefore there exist $\varepsilon_0$ and $\varepsilon_1$ s.t. the
matrix $\p_\mu G(\mu,v)$ has a uniformly bounded inverse for
$\mu\in[\frac{1}{4},1]\times[0,\varepsilon_0]$ and  $v\in
\bigcup_{\mu\in[\frac{1}{4},1]\times[0,\varepsilon_0]}
B_{\varepsilon_1}(V_{\mu})$. Hence by the implicit function theorem,
the equation $G(\mu,v)=0$ has a unique solution $\mu=g(v)$ on a
neighborhood of every $V_\mu$,
$\mu\in[\frac{1}{4},1]\times[0,\varepsilon_0]$, which is $C^1$ in
$v$.  Our next goal is to determine these neighborhoods.

To determine a domain of the function $\mu=g(v)$, we examine closely
a proof of the implicit function theorem. Proceeding in a standard
way, we expand the function $G(\mu,v)$ in $\mu$ around $\mu_0$:
\begin{equation*}
G(\mu,v)=G(\mu_0,v)+\p_\mu G(\mu_0,v)(\mu-\mu_0)+R(\mu,v),
\end{equation*}
where $R(\mu,v)=\O{|\mu-\mu_0|^2}$ uniformly in $v\in X$. Here
$|\mu|^2=|a|^2+|b|^2$ for $\mu=(a,b)$.  Inserting this into the
equation $G(\mu,v)=0$ and inverting the matrix $\p_\mu G(\mu_0,v)$,
we arrive at the fixed point problem $\alpha=\Phi_v(\alpha)$, where
$\alpha:=\mu-\mu_0$ and $\Phi_v(\alpha):=- \p_\mu G(\mu_0,v)^{-1}
[G(\mu_0,v)+R(\mu,v)]$.  By the above estimates there exists an
$\varepsilon_1$ such that the matrix $\p_\mu G(\mu_0,v)^{-1}$ is
bounded uniformly in $v\in B_{\varepsilon_1}(V_{\mu_0})$.  Hence we
obtain from the remainder estimate above that
\begin{equation}
|\Phi_v(\alpha)|\lesssim|G(\mu_0,v)|+|\alpha|^2.
\label{eqn:SplittingSharp}
\end{equation}
Furthermore, using that $\p_\alpha \Phi_v(\alpha)= - \p_\mu
G(\mu_0,v)^{-1} [G(\mu,v)- G(\mu_0,v)+R(\mu,v)]$ we obtain that there exist
$\varepsilon \leq \varepsilon_1$ and $\delta$ such that
$\|\p_\alpha\Phi_v(\alpha)\|\le\frac{1}{2}$ for all $v\in
B_{\varepsilon}(V_{\mu_0})$ and $\alpha\in B_\delta(0)$.  Pick
$\varepsilon$ and $\delta$ so that $\varepsilon
\ll\delta\ll b_0\ll 1$.  Then, for all $v\in
B_{\varepsilon}(V_{\mu_0})$, $\Phi_v$ is a contraction on the ball
$B_\delta(0)$ and consequently has a unique fixed point in this
ball.  This gives a $C^1$ function $\mu=g(v)$ on
$B_{\varepsilon}(V_{\mu_0})$ satisfying $|\mu-\mu_0|\le\delta$. An
important point here is that since $\varepsilon\ll b(0)$ we have
that $b>0$ for all $V_{ab}\in B_{\varepsilon}(V_{\mu_0})$.  Now,
clearly, the balls $B_{\varepsilon}(V_{\mu_0})$ with
$\mu_0\in[\frac{1}{4},1]\times[0,\varepsilon_0]$ cover the
neighbourhood $U_{\varepsilon_0}$.  Hence, the map $g$ is defined on
$U_{\varepsilon_0}$ and is unique, which implies the first part of
the proposition.

Now we prove the second part of the proposition.  The definition of
the function $G(\mu,v)$ implies $G(\mu_0,v)=G(\mu_0,v-V_{\mu_0})$
and
\begin{equation}
|G(\mu_0,v)|\lesssim \|e^{-\frac{1}{9} y^2}(v-V_{\mu_0})\|_\infty.
\label{eqn:28aA}
\end{equation}
This inequality together with the estimate
\eqref{eqn:SplittingSharp} and the fixed point equation
$\alpha=\Phi_v(\alpha)$, where $\alpha=\mu-\mu_0$ and $\mu=g(v)$,
implies $|\alpha|\lesssim \|e^{-\frac{1}{9}
y^2}(v-V_{\mu_0})\|_\infty+|\alpha|^2$ which, in turn, yields
\eqref{eqn:28a}.
\end{proof}

\begin{prop}
\label{Prop:SplittingIC} In the notation of Proposition
\ref{Prop:Splitting}, if $\|\langle y\rangle^{-n}(v-V_{a_0
b_0})\|_\infty\le \delta_{n}$ with $n=0,3$,
$\delta_{3}=O(b^{2}_{0})$ and $\delta_{0}$ small, then
\begin{equation}\label{eq:vv0}
|g(v)-(a_0, b_0)|\lesssim b_{0}^{2},
\end{equation}
\begin{equation}
\|\langle y\rangle^{-3}(v-V_{g(v)}))\|_\infty\lesssim b_{0}^{2}
\label{Ineq:IC}
\end{equation}
and
\begin{equation}\label{eq:vwithoutweight}
\|v-V_{g(v)}\|_\infty\lesssim \delta_{0}+b_{0}.
\end{equation}
\end{prop}
\begin{proof}
Let $g(v)=(a,b)$ and $\mu=(a_0,b_0)$.   By
\eqref{eqn:SplittingSharp} and the fixed point equation
$\alpha=\Phi_v(\alpha)$, we have $\alpha\lesssim
|G(\mu_0,v)|+|\alpha|^2$ which, in turn, yields $|\mu-\mu_0|\lesssim
|G(\mu_0,v)|$.  By \eqref{eqn:28aA} and one of the conditions of the
proposition, $G(\mu_0,v)=O(b_0^2)$ if $a_0\in[\frac{1}{4},1]$. The
last two estimates imply \eqref{eq:vv0}.  Using Equation
\eqref{eqn:28a} we obtain
$$
\begin{array}{lll}
\|\langle y\rangle^{-3}(v-V_{g(v)})\|_\infty&\leq &\|\langle
y\rangle^{-3}(v-V_{\mu_{0}})\|_\infty+\|\langle
y\rangle^{-3}(V_{g(v)}-V_{\mu_{0}})\|_\infty\\
&\lesssim& \|\langle
y\rangle^{-3}(v-V_{\mu_{0}})\|_\infty+|g(v)-\mu_{0}|\\
&\lesssim & \|\langle y\rangle^{-3}(v-V_{\mu_{0}})\|_\infty
\end{array}
$$ which leads to (~\ref{Ineq:IC}).
Finally, to prove Equation (~\ref{eq:vwithoutweight}),
we write
$$\|v-V_{g(v)}\|_\infty\leq \|v-V_{a_{0},b_{0}}\|_\infty+\|V_{g(v)}-V_{a_{0},b_{0}}\|_\infty.$$
A straightforward computation gives $\|V_{a b}-V_{a_{0}
b_{0}}\|_\infty\lesssim |a-a_{0}|+\frac{|b-b_{0}|}{b_{0}}$.  Since
by (~\ref{eq:vv0}), $|a-a_{0}|+|b-b_{0}|=O(b_{0}^2)$, we have
$\|V_{a b}-V_{a_0 b_0}\|_\infty\lesssim b_{0}$. This together with
the fact $\|v-V_{a_{0},b_{0}}\|_\infty\leq \delta_{0}$ completes the
proof of (~\ref{eq:vwithoutweight}).
\end{proof}
\section{A priori Estimates}
\label{Section:APriori} In this section we assume that (~\ref{NLH})
has a unique solution, $u(x,t),$ $0\leq t\leq t_{*}$, such that
$v(y,\tau)=\lambda^{-\frac{2}{p-1}}(t)u(x,t)$, where $y=\lambda x$
and $\tau(t):= \int_{0}^{t}\lambda^{2}(s)ds$, is in the neighborhood
$U_{\epsilon_{0}}$ determined in Proposition ~\ref{Prop:Splitting}.
Then by Proposition ~\ref{Prop:Splitting} there exist $C^{1}$
functions $a(\tau)$ and $b(\tau)$ such that $v(y,\tau)$ can be
represented as
\begin{equation}\label{eqn:split2}
v(y,\tau)=\lb\frac{2c}{p-1+by^{2}}\rb^{\frac{1}{p-1}}+e^{\frac{ay^{2}}{4}}\xi(y,\tau)
\end{equation}
where $\xi(\cdot,\tau)\perp \phi_{0a},\phi_{2a}$ (see
(~\ref{eqn:split})) and $c=\frac{1}{2}a+\frac{1}{4}$ and, by
Proposition ~\ref{Prop:SplittingIC},
\begin{equation}\label{eqn:splitB}
\|\lra{y}^{-3}e^{\frac{ay^{2}}{4}}\xi\|_\infty\lesssim b^2(\tau).
\end{equation} Now we set
$$\lambda^{-3}(t)\partial_{t}\lambda(t)=a(\tau(t)).$$ In this
section we present a priori bounds on the fluctuation $\xi$ which
are proved in later sections.

We begin with defining convenient estimating functions. Denote by
$\chi_{\geq D}$ and $\chi_{\leq D}$ the characteristic functions of
the sets $\{|x|\geq D\}$ and $\{|x|\leq D\}:$
\begin{equation}\label{cutoff} \chi_{\geq D}(x):=
\left\{\begin{array}{lll}
1\ \text{if}\ |x|\geq D\\
0\ \text{otherwise}
\end{array}\right.
\ \text{and}\ \chi_{\leq D}:=1-\chi_{\geq D}.
\end{equation} We take $D:=\frac{C}{\sqrt{\beta}}$ where
$C$ is a large constant to be specified in Section \ref{SEC:EstM2}.
Let the function $\beta(\tau)$ and the constant $\kappa$ be defined
as
\begin{equation}\label{FunBTau}
\beta(\tau):=\frac{1}{\frac{1}{b(0)}+\frac{4p}{(p-1)^{2}}\tau}\
\mbox{and}\ \kappa:=\min\{\frac{1}{2},\frac{p-1}{2}\}.
\end{equation}
For the functions $\xi(\tau),$ $b(\tau)$ and $a(\tau)$ we introduce
the following estimating functions (families of semi-norms)
\begin{equation}
\label{majorants}
\begin{array}{lll}
M_{1}(T)&:=\max_{\tau\leq T}
\beta^{-2}(\tau)\|\lra{y}^{-3}e^{\frac{a}{4}y^2}\xi(\tau)\|_{\infty},\\
M_{2}(T)&:=\max_{\tau\leq
T}\|e^{\frac{a}{4}y^{2}}\chi_{\geq D}\xi(\tau)\|_{\infty},\\
A(T)&:=\max_{\tau\leq
T}\beta^{-2}(\tau)\left|a(\tau)-\frac{1}{2}+\frac{2b(\tau)}{p-1}\right|,\\
B(T)&:=\max_{\tau\leq
T}\beta^{-(1+\kappa)}(\tau)|b(\tau)-\beta(\tau)|.
\end{array}
\end{equation}
\begin{prop}\label{Prop:aprior}
Let $\xi$ to be defined in (~\ref{eqn:split2}) and assume
$M_{1}(0),A(0), B(0)\lesssim 1$, $M_{2}(0)\lll 0$. Assume there
exists an interval $[0,T]$ such that for $\tau\in [0,T]$,
$$
M_{1}(\tau), A(\tau),\ B(\tau)\leq \beta^{-\kappa/2}(\tau).$$ Then
in the same time interval the parameters $a$, $b$ and the function
$\xi$ satisfy the following estimates
\begin{equation}\label{EstB}
|b_{\tau}(\tau)+\frac{4p}{(p-1)^{2}}b^{2}(\tau)| \lesssim
\beta^{3}(\tau)+\beta^{3}(\tau)M_{1}(\tau)(1+A(\tau))+\beta^{4}(\tau)M_{1}^{2}(\tau)+\beta^{2p}M_{1}^{2p}(\tau),
\end{equation}
and
\begin{equation}\label{MajorE}
B(\tau)\lesssim
1+M_{1}(\tau)(1+A(\tau))+M_{1}^{2}(\tau)+M_{1}^{p}(\tau),
\end{equation}
\begin{equation}\label{EstA}
A(\tau)\lesssim
A(0)+1+\beta(0)M_1(\tau)(1+A(\tau))+\beta(0)M_{1}^{2}(\tau)+\beta^{2p-2}(0)M_{1}^{p}(\tau),
\end{equation}
\begin{equation}\label{M1}
\begin{array}{lll}
M_{1}(\tau)
&\lesssim &  M_{1}(0)+\beta^{\frac{\kappa}{2}}(0)[1+M_{1}(\tau)A(\tau)+M_{1}^{2}(\tau)+M_{1}^{p}(\tau)]\\
& &+[M_{2}(\tau)M_{1}(\tau)+M_{1}(\tau)M_{2}^{p-1}(\tau)],
\end{array}
\end{equation}
\begin{equation}\label{M2}
\begin{array}{lll}
M_{2}(\tau)&\lesssim& M_{2}(0)+\beta^{1/2}(0)M_{1}(0)+M_{2}^{2}(\tau)+M_{2}^{p}(\tau)\\
&
&+\beta^{\frac{\kappa}{2}}(0)[1+M_{2}(\tau)+M_{1}(\tau)A(\tau)+M_{1}^{2}(\tau)+M_{1}^{p}(\tau)].
\end{array}
\end{equation}
\end{prop}

Equations \eqref{EstB}-\eqref{EstA} will be proved in Section
\ref{SEC:EstB}.  Equations \eqref{M1} and \eqref{M2} will be proved
in Sections ~\ref{SEC:EstM1} and ~\ref{SEC:EstM2} respectively.

\section{Proof of Main Theorem \ref{maintheorem}}\label{SecMain}
We begin with an analysis of the initial conditions. In the next
lemma we show that restriction (~\ref{INI2}) on the initial
conditions involving two parameters can be rescaled into a condition
involving one parameter.
\begin{lemma}\label{LM:rescale}
Let $u_{0}$ satisfy the condition (~\ref{INI2}) and let
$k_{0}:=(\sqrt{2\cO+\frac{2}{p-1}b_{0}})^{-1/2}$ and
$\beta_{0}:=b_{0}k_{0}^{2}$. Then we have the estimates
$$
\|\langle
k_{0}x\rangle^{-n}(k_{0}^{\frac{2}{p-1}}u_{0}(k_{0}x)-(\frac{1-\frac{2}{p-1}\beta_{0}}{p-1+\beta_{0}
x^{2}})^{\frac{1}{p-1}})\|_{\infty}\leq \delta_{n},\ n=0,3;\
\delta_3 =  C\beta_{0}^2.$$
\end{lemma}
\begin{proof}
It is straightforward to verify that the function
$k_{0}^{\frac{2}{p-1}}u_{0}(k_{0}x)$ has all the properties above.
\end{proof}
Due to this lemma, in what follows, it suffices to consider the
condition \eqref{INI2} of Theorem \ref{maintheorem} with $c_0=
\frac{1}{2}-\frac{1}{p-1}b_{0}$ (since $2c = a + \frac{1}{2}$, this
gives $a_0= \frac{1}{2}-\frac{2}{p-1}b_{0}$):
\begin{equation}\label{INI}
\|\langle
x\rangle^{-n}(u_{0}(x)-(\frac{1-\frac{2}{p-1}b_{0}}{p-1+b_{0}x^{2}})^{\frac{1}{p-1}})\|_{\infty}\le
\delta_{n},\ n=0,3,\ \delta_3 = C b_{0}^2.
\end{equation}
To obtain the statement of Theorem \ref{maintheorem} we rescale the
result obtained below as
\begin{equation}\label{resc}
u_{thm}(x, t)= k_{0}^{-\frac{2}{p-1}}u_{pf}(x/k_{0}, t/k_{0}^2).
\end{equation}
Here $u_{thm}(x, t)$ and $u_{pf}(x, t)$ are the solutions $u(x, t)$
appearing in the theorem and in the proof, respectively. This
rescaling and the constrain $c=\frac{1}{2}a + \frac{1}{4}$ used in
the proof give the following relations between the parameters used
in the theorem and the proof:
\begin{equation}\label{param}
\lambda_{thm}( t)= k_{0}^{-1}\lambda_{pf}( t),\  b_{thm}( t)=
b_{pf}(\tau( t)),\ c_{thm}( t)= c_{pf}(\tau( t)).
\end{equation}

By Theorem ~\ref{THM:Local} (the local existence theorem, proven in
Appendix \ref{Appendix:LWP}), there exists $\infty\ge t_{*}>0$ such
that Equation (~\ref{NLH}) has a unique solution $u(x,t)$ in
$C([0,t_{*}],L^{\infty})$ and, if $t_* < \infty$, then
$\|u(\cdot,t)\|_{\infty} \rightarrow \infty$ as $t\rightarrow t_*$.
Recall that the solutions $u(x,t)$, $v(y, \tau)$ and $w(y,\tau)$ and
the corresponding initial conditions are related by the scaling and
gauge transformations (see (~\ref{eq:definev}) and
(~\ref{eq:definew})). Take $\lambda(0)=1$. Then we have that
$u_{0}(x)=v_{0}(y)$.

Choose $b_0$ so that $Cb_{0}^2\leq \frac{1}{2}\epsilon_{0}$ with $C$
the same as in \eqref{INI} and with $\epsilon_{0}$ given in
Proposition ~\ref{Prop:Splitting}. Then $v_{0}\in
U_{\frac{1}{2}\epsilon_{0}}$, by the condition \eqref{INI} on the
initial conditions with $n=3$.
By continuity there is a (maximal) time $t_\#\leq t_{*}$ such that
$v \in U_{\epsilon_0}$ for $t< t_\#$. For this time interval
Propositions \ref{Prop:Splitting} and \ref{Prop:SplittingIC} hold
for $v$ and, in particular, we have the splitting
\eqref{eqn:split2}. Recall that we assume $a=2c-\frac{1}{2}$ in the
decomposition \eqref{eqn:split2}. This implies that the initial
condition can be written in the form
\begin{equation} \label{eqn:splitic}
v_0(y) =V_{a(0) b(0)}(y) + e^{\frac{a(0)y^{2}}{4}}\xi_0(y),
\end{equation}
where $(a(0), b(0))=g(v_0)$ and $\xi_0\perp e^{-\frac{a(0)}{4} y^2},
(1-a(0) y^2)e^{-\frac{a(0)}{4} y^2}$.

By the relation $\beta(0)=b(0)$, Equation (~\ref{INI}) and
Proposition ~\ref{Prop:SplittingIC},   $A(0)$, $M_{1}(0)\lesssim 1$
and $M_{2}(0) \ll 1$,
while $B(0)=0$, by the definition. Since $\beta(\tau) \le \beta(0)
\ll 1$,  we have, by the continuity (or by Proposition
~\ref{Prop:SplittingIC}), that for a sufficiently small time
interval
\begin{equation}\label{ApriorEST}
M_{1}(\tau), \ B(\tau),\ A(\tau)\leq
\beta^{-\frac{\kappa}{2(2+p)}}(0)\leq
\beta^{-\frac{\kappa}{2}}(\tau),
\end{equation}
where, recall, the definitions of $\beta(\tau)$ and $\kappa$ are
given in (~\ref{FunBTau}). Then Equations
(~\ref{MajorE})-(~\ref{M2}) imply that for the same time interval
\begin{equation}\label{UtimateEST}
M_{1}(\tau),\ B(\tau),\ A(\tau)\lesssim 1, \ M_{2}(\tau)\ll 1.
\end{equation}
(In fact, $M_{i}(\tau) \lesssim M_{i}(0)  +
\beta^{\frac{\kappa}{2}}(0),\ i=1,2$.) Indeed, since
$M_{1}(\tau)\leq \beta^{-\frac{\kappa}{2}}(0)$, we can solve
\eqref{EstA} for $A(\tau)$. We substitute the result into Eqns
\eqref{M1} - \eqref{M2} to obtain inequalities involving only the
estimating functions $M_{1}(\tau)$ and $M_{2}(\tau)$. Consider the
resulting inequality for $M_{2}(\tau)$. The only terms on the
r.h.s., which do not contain $\beta(0)$ to a power at least
$\kappa/2$ as a factor, are $M_{2}^{2}(\tau)$ and $M_{2}^{p}(\tau)$.
Hence for $M_{2}(0)\ll 1$ this inequality implies that $M_{2}(\tau)
\lesssim M_{2}(0)  + \beta^{\frac{\kappa}{2}}(0)$. Substituting this
result into the inequality for $M_{1}(\tau)$ we obtain that
$M_{1}(\tau) \lesssim M_{1}(0)  + \beta^{\frac{\kappa}{2}}(0)$ as
well. The last two inequalities together with \eqref{MajorE} and
\eqref{EstA} imply the desired estimates on $A(\tau)$ and $B(\tau)$.

By \eqref{UtimateEST} and continuity, \eqref{ApriorEST} holds on a
larger interval which in turn implies \eqref{UtimateEST} on this
larger time interval and so forth.  Hence, \eqref{UtimateEST} holds
for $t<t_\#=t_*$.

By the definitions of $A(\tau)$ and $B(\tau)$ in (~\ref{majorants})
and the facts that $A(\tau), B(\tau)\lesssim 1 $ proved above and
the relation $2c=a+\frac{1}{2}$, we have that
\begin{equation}\label{EstABtau}
a(\tau)-\frac{1}{2}=-\frac{2}{p-1}b(\tau)+\O{\beta^{2}(\tau)},\
b(\tau)=\beta(\tau)+\O{\beta^{1+\kappa/2}(\tau)}.
\end{equation}
Hence $a(\tau)-\frac{1}{2}=O(\beta(\tau))$.  Recall that
$a=\lambda^{-3}\frac{d}{dt}\lambda$, which can be rewritten as
$\lambda(t)^{-2}=1-2 \int_{0}^{t}a(\tau(s))ds$ or
\begin{equation}
\lambda(t)=[1-2 \int_{0}^{t}a(\tau(s))ds]^{-1/2}
\end{equation} where we use $\lambda(0)=1.$ Since $|a(\tau(t))-\frac{1}{2}|=O(b(\tau(t)))$, there exists a time $t^{*}$ such
that $1=2 \int_{0}^{t^*}a(\tau(s))ds$, i.e. $\lambda(t)\rightarrow
\infty$ as $t\rightarrow t^{*}$. Furthermore, by the definition of
$\tau$ and the property of $a$ we have that $\tau(t)\rightarrow
\infty$ as $t\rightarrow t^{*}.$ We will show below that $t^{*}$ is
the blow-up time.

Equation (~\ref{EstABtau}) implies $b(\tau(t))\rightarrow 0$ and
$a(\tau(t))\rightarrow \frac{1}{2}$ as $t\rightarrow t^*$.  By the
analysis above and the definitions of $a$, $\tau$ and $\beta$ (see
\eqref{FunBTau}) we have
$$\lambda(t)=(t^*-t)^{-1/2}(1+o(1)),\ \tau(t)=-ln|t^*-t|(1+o(1)),$$
and $$\beta(\tau(t))=-\frac{(p-1)^{2}}{4pln|t^*-t|}(1+o(1)).$$ By
(~\ref{EstABtau}) we have
$$b(\tau(t))=\frac{(p-1)^{2}}{4p(|ln|t^*-t||)}(1+O(|\frac{1}{ln|t^*-t|}|^{\kappa/2}))$$
and
$$a(\tau(t))=\frac{1}{2}-\frac{p-1}{2p|ln|t^*-t||}(1+O(\frac{1}{ln|t^*-t|})).$$
The last equation together with the relation $c=\frac{1}{2}a +
\frac{1}{4}$ implies
$$c(\tau(t))=\frac{1}{2}-\frac{p-1}{4p|ln|t^*-t||}(1+O(\frac{1}{ln|t^*-t|})).$$

Now, using the relation between the functions  $u(x,t)$ and
$v(y,\tau)$ and the splitting result (Proposition
~\ref{Prop:Splitting}) we obtain the following a priory estimate on
the (non-rescaled) solution $u(x,t)$ of equation (~\ref{NLH}):
$$ \|u(t)\|_{\infty} \lesssim \lambda(t)^{\frac{2}{p-1}} [1 + M_1(\tau)
+M_2(\tau)],$$ where $\tau=\tau(t)$ is defined above and we use the
fact $\|e^{\frac{ay^{2}}{4}}\xi(\cdot,\tau)\|_{\infty}\lesssim
M_{1}(\tau)+M_{2}(\tau)$. By the estimate (~\ref{UtimateEST}) above
the majorants $M_{j}(\tau)$ are uniformly bounded and therefore
\begin{equation}\label{upperbound}
\|u(t)\|_{\infty} \lesssim  \lambda(t)^{\frac{2}{p-1}}
\end{equation} for $t\leq t_{*}.$

Recall that if $t_* < \infty$, then $\|u(\cdot,t)\|_{\infty}
\rightarrow \infty$ as $t\rightarrow t_*$. Hence, $t_*\ge t^*$, by
the bounds \eqref{upperbound} and $\lambda(t) < \infty$ for $t <
t^*$. On the other hand, if $t_*>t^*$, then by \eqref{eqn:split2}
and \eqref{eqn:splitB}
\begin{equation}
|u(0,t)|\ge \lambda(\tau(t))^\frac{2}{p-1}\lsb \lb\frac{2
c(\tau(t))}{p-1}\rb^\frac{1}{p-1}-C b(\tau(t))^2
\rsb\rightarrow\infty,
 \label{eqn:sqeefdsfew}
\end{equation}
as $t\uparrow t^*$, which contradicts the existence of $u(x,t)$ on
$[0,t_*)$.  Hence $t_{*}=t^{*}$. Thus we have shown the existence of
the solution $u$ up to the time $t^{*}$ having $v \in
U_{\epsilon_0}$ and obeying the estimates \eqref{UtimateEST}. Then
the results above describing the dynamics of the parameters $a$,
$b$, $c$ and $\lambda$ as well as \eqref{UtimateEST} imply that $u$
blows up at the time $t^{*}$, see Equation \eqref{eqn:sqeefdsfew}.
Furthermore, Equations \eqref{eqn:split2} and \eqref{UtimateEST}
imply, after rescaling \eqref{resc}, the second  statement of
Theorem ~\ref{maintheorem}. Finally, the third statement follows
from the asymptotic expressions for the parameters $b$, $c$ and
$\lambda$ obtained above and the relations in \eqref{param}. This
completes the proof of Theorem ~\ref{maintheorem}.
\begin{flushright}
$\square$
\end{flushright}
\section{Lyapunov-Schmidt Splitting (Effective Equations)}
\label{Section:Splitting}
According to Proposition \ref{Prop:Splitting} the solution
$w(y,\tau)$ of \eqref{eqn:w} can be decomposed as
\eqref{eqn:split2}, with the parameters $a$, $b$ and $c$ and the
fluctuation $\xi$ depending on time $\tau$:
\begin{equation}
w=v_{abc}+\xi,\ \xi\bot\ \phi_{0,a},\ \phi_{2 a}, \label{eqn:split3}
\end{equation}
where, recall, $v_{a b c}:=v_{c b } e^{-\frac{a}{4} y^2}$ and
$c=\frac{1}{2}a+\frac{1}{4}$.  In this section we derive equations
for the parameters $a(\tau)$, $b(\tau)$ and $c(\tau)$ and the
fluctuation $\xi(y,\tau)$.

Plugging the decomposition \eqref{eqn:split3} into \eqref{eqn:w}
gives the equation
\begin{equation} \xi_\tau=-{\mathcal L}_{abc}\xi+{\mathcal
N}(\xi,a,b,c) +{\mathcal F}(a,b,c), \label{eqn:FluctuationGauged}
\end{equation}
where the operator ${\mathcal L}_{a,b,c}$, the functions
$\mathcal{N}(\xi,a,b,c)$ and ${\mathcal F}(a,b,c)$ are defined as
\begin{align}
{\mathcal L}_{a b c}:=&-\p_y^2+\frac{1}{4} \lb a^2 + a_\tau\rb
y^2-\frac{a}{2}+\frac{2 a}{p-1}-\frac{2 p c}{p-1+ b y^2},\label{eqn:LinOp}\\
{\mathcal N}(\xi,b,c):=
&\left[|\xi+v_{a b c}|^{p-1}(\xi+v_{a b c})-v_{a b c}^p-p v_{a b c
}^{p-1}\xi\right]e^{\frac{a}{4}(p-1)y^2},\\
{\mathcal F}(a,b,c):=&\frac{1}{p-1}\left[
\Gamma_{0}+\Gamma_{1}\frac{(p-1)ay^{2}}{p-1+by^{2}}-\frac{4pb^{3}y^{4}}{(p-1)^{2}(p-1+by^{2})^{2}}
\right] v_{a b c},\label{eqn:Forcing}
\end{align}
with the functions $\Gamma_{0}$ and $\Gamma_{1}$ given as
\begin{align}
\Gamma_0&:=-\frac{c_\tau}{c}+2(c-a)-\frac{2}{p-1} b,\label{eqn:Gamma0}\\
\Gamma_1&:=\frac{1}{a(p-1)}\lb b_\tau-2
b(c-a)+\frac{2(3p-1)}{(p-1)^2}b^2 \rb.\label{eqn:Gamma2}
\end{align}
\begin{prop}
\label{Prop:ForcingBounds} If $A(\tau),\ B(\tau)\leq
\beta^{-\frac{\kappa}{2}}(\tau)$ and $1/4\leq c(0)\leq 1$, then
\begin{align}
\|\lra{y}^{-3} e^{\frac{a}{4} y^2}{\cal
F}\|_\infty=\O{|\Gamma_0|+|\Gamma_1|+\beta^\frac{5}{2}}\ \mbox{and}\
\|e^{\frac{a}{4} y^2}{\cal
F}\|_\infty=\O{|\Gamma_0|+\frac{1}{\beta}|\Gamma_1|+\beta}.
\label{eqn:Festimate}
\end{align}
Furthermore we have for ${\cal N}={\cal N}(\xi,b,c)$
\begin{equation}
|{\cal N}|\lesssim
e^{\frac{ay^{2}}{4}}|\xi|^{2}+e^{(p-1)\frac{a}{4}y^{2}}|\xi|^{p}.
\label{eqn:69a}
\end{equation}
\end{prop}
\begin{proof}
Rearranging the leading term of expression for $\mathcal{F}$ so that
$y^2$ appears in the combination $a y^2-1$ gives the more convenient
expression
\begin{equation}
{\cal F}=\frac{1}{p-1}\lsb\Gamma_0+\Gamma_1+\Gamma_1 (a
y^2-1)-\Gamma_{1}\frac{aby^{4}}{p-1+by^{2}}+G_{1}\rsb v_{a b
 c}\label{eqn:ForcingDecomposed}
\end{equation} with $G_{1}:=-\frac{4pb^{3}y^{4}}{(p-1)^{2}(p-1+by^{2})^{2}}.$ We estimate $\|\lra{y}^{-3} e^{\frac{a}{4}
y^2}{\cal F}\|_\infty$ using this form of ${\cal F}$ and the
estimates
\begin{align*}
\| e^{\frac{a}{4} y^2} v_{a b c}\|_\infty,\ \|\lra{y}^{-3} (a y^2-1)
\|_\infty,\ \lesssim 1.
\end{align*}  The result is
\begin{equation}
\|\lra{y}^{-3} e^{\frac{a}{4} y^2}{\cal F}\|_\infty\lesssim
|\Gamma_0|+(1+ b^\frac{1}{2})|\Gamma_1|+b^\frac{5}{2}.
\label{eqn:singleast}
\end{equation}

The estimate of $\|e^{\frac{a}{4}y^2}{\cal F}\|$ is proved in a
similar way as the first estimate. Recall the expression of
$\mathcal{F}$ in Equation \eqref{eqn:Forcing}. We use the estimates
\begin{align*}
\left\|e^{\frac{a}{4} y^2} v_{a b c}\right\|_\infty,\
 \left\|e^{\frac{a}{4} y^2}\frac{(b y^2)^n}{(p-1+b
y^2)^2} v_{a b c}\right\|_\infty&\lesssim 1,\ n=0,1,
\end{align*}
to obtain that
\begin{equation}
\|e^{\frac{a}{4} y^2}{\cal F}\|_\infty\lesssim
|\Gamma_0|+\frac{1}{b}|\Gamma_1|+b . \label{eqn:doubleast}
\end{equation}

Now we estimate $b$ in terms of $\beta$ and $B$ to complete the
proof of the first bound. The assumption that $B\leq
\beta^{-\frac{\kappa}{2}}$ implies that $b=
\beta+\O{\beta^{1+\frac{\kappa}{2}}}$, which together with estimates
\eqref{eqn:singleast} and \eqref{eqn:doubleast}, implies the first
estimate \eqref{eqn:Festimate}.

For (~\ref{eqn:69a}) we observe that if $v_{abc}\le 2|\xi|$ then
$|{\cal N}|\le e^{(p-1)\frac{a}{4} y^2}(p+3)|\xi|^p$.  If
$v_{abc}\ge 2|\xi|$, then we use the formula ${\cal
N}=e^{(p-1)\frac{a}{4} y^2}p \int_0^1 \lsb(v_{a b c}+s\xi)^{p-1}-v_{a
b c}^{p-1}\rsb\xi\, ds$ and consider the cases $1<p\le 2$ and $p>2$
separately to obtain \eqref{eqn:69a}.
\end{proof}


\begin{prop}\label{Prop:LyapunovSchmidtReduction}
Recall that $a=2c-\frac{1}{2}$.  Suppose that $A(\tau), B(\tau),
M_{1}(\tau)\leq \beta^{-\frac{\kappa}{2}}$ and $1/4\leq c(0)\leq 1$
for $0\le\tau\le T$. Let $w=v_{a b c}+\xi$ be a solution to
\eqref{eqn:w} with $\xi\bot \phi_{0a},\ \phi_{2a}$. Over times
$0\le\tau\le T$, the parameters $b$ and $c$ satisfy
\begin{align}
b_\tau&=-\frac{2(3p-1)}{(p-1)^2}b^2+2 b(c-a)+{\cal R}_b(\xi,b,c),\label{eqn:bParameter}\\
\frac{c_\tau}{c}&=2(c-a)-\frac{2}{p-1} b+{\cal
R}_c(\xi,b,c),\label{eqn:cParameter}
\end{align}
where the remainders ${\cal R}_b$ and ${\cal R}_c$ are of the order
$\O{\beta^3+\beta^3 M_1(1+
A)+\beta^{4}M_{1}^{2}+\beta^{2p}M_{1}^{p}}$ and satisfy ${\cal
R}_{b}(0,b,c), {\cal R}_{c}(0,b,c)=O(b^3).$
\end{prop}
\begin{proof} We take inner product of the equation
(~\ref{eqn:FluctuationGauged}) with $\phi_{ja}$ to get
$$\langle \xi_{\tau}, \phi_{ja}\rangle=\langle -{\mathcal L}_{abc}\xi+{\mathcal
N}(\xi,a,b,c) +{\mathcal F}(a,b,c),\phi_{ja}\rangle.$$  We use the
orthogonality conditions $\phi_{ja}\perp \xi$ to derive
(~\ref{eqn:bParameter}) and (~\ref{eqn:cParameter}).  We start with
analyzing the $\mathcal{F}$ term. The inner product of
\eqref{eqn:ForcingDecomposed} with $\phi_{0a}$ and $\phi_{2a}$ gives
the expression
\begin{multline}
(p-1)\ip{{\cal F}}{\phi_{ja}}=(\Gamma_0+\Gamma_1)\ip{v_{a b
c}}{\phi_{ja}}+\Gamma_1\ip{v_{a b c}}{(a
y^2-1)\phi_{ja}}-\Gamma_{1}\langle \frac{aby^{4}}{p-1+by^{2}}v_{a b
c},\phi_{ja}\rangle+\langle G_{1}v_{a b c},\phi_{ja}\rangle
\label{eqn:ProjectedEquation}
\end{multline}
where $j=0$ or $2$.  By rescaling the variable of integration so
that the exponential term does not contain the parameter $a$,
expanding $v_{a b c}$ to the constant term in $\frac{b}{a}$ and
estimating the remainder by $\O{a^{-\frac{1}{2}} b y^2 e^{-y^2/2}}$
we obtain the estimates
\begin{align*}
\ip{v_{a b
c}}{\phi_{0a}}&=(\frac{2c}{p-1})^{\frac{1}{p-1}}\sqrt{\frac{2\pi}{a}}+\O{b},\\
\ip{v_{a b c}}{\phi_{2a}}&=\O{b},\\
\ip{v_{a b c}}{(a y^2-1)
\phi_{2a}}&=(\frac{2c}{p-1})^{\frac{1}{p-1}}\sqrt{\frac{8\pi}{a}}+\O{b},
\end{align*}
$$
\ip{v_{a b c}}{\lra{y}^3 \phi_{0a}},\ \ip{v_{a b c}}{\lra{y}^3
\phi_{2a}}\lesssim 1.
$$
Substituting these estimates into Equations
\eqref{eqn:ProjectedEquation} and recalling the definition of $G_1$
gives
\begin{align}
\ip{{\cal F}}{\phi_{0a}}&=\frac{1}{p-1}(\frac{2c}{p-1})^{\frac{1}{p-1}}\sqrt{\frac{2\pi}{a}}(\Gamma_0+\Gamma_1)+R_1,\label{est:ForcingPsi0}\\
\ip{{\cal
F}}{\phi_{2a}}&=\frac{1}{p-1}(\frac{2c}{p-1})^{\frac{1}{p-1}}\sqrt{\frac{8\pi}{a}}\Gamma_1+R_2,\label{est:ForcingPsi2}
\end{align}
where both remainders $R_1$ and $R_2$ are bounded by $\O{b
|\Gamma_0|+b |\Gamma_1|+b^3}$.

To estimate the projection of $\p_\tau\xi$ onto $\phi_{0a}$ and
$\phi_{2a}$, we differentiate the orthogonality conditions
$\ip{\xi}{\phi_{0a}}=0$ and $\ip{\xi}{\phi_{2a}}=0$, obtaining the
relations $\ip{\xi_\tau}{\phi_{0a}}=-\ip{\xi}{\p_\tau\phi_{0a}}$ and
$\ip{\xi_\tau}{\phi_{2a}}=-\ip{\xi}{\p_\tau\phi_{2a}}$.  When
simplified using the orthogonality conditions on $\xi$, these
relations give
\begin{equation*}
\ip{\xi_\tau}{\phi_{0a}}=0\ \mbox{and}\
|\ip{\xi_\tau}{\phi_{2a}}|\leq
|\frac{1}{4}a^{-1}a_\tau\ip{\lra{y}^{-3} e^{\frac{a}{4}
y^2}\xi}{a^2\lra{y}^{3} y^4 e^{-\frac{a}{2} y^2}}|.
\end{equation*}
Estimating the right hand side of the second inequality by
H\"{o}lder's inequality and using the definition of $M_1(\tau)$
gives that over times $0\le \tau\le T$
\begin{equation*}
\ip{\xi_\tau}{\phi_{2a}}=\O{|a_\tau| \beta^2 M_1}.
\end{equation*}
Next we replace $a_\tau$ in with expressions involving $\Gamma_0$
and $\Gamma_1$.  Since $a=2c-\frac{1}{2}$, $a_\tau=2c_{\tau}$. From
\eqref {eqn:Gamma0} and \eqref{eqn:Gamma2},
$$
c_\tau=\O{\Gamma_0+\beta^2 A}
$$
for times $0\le \tau\le T$.  Substituting these estimates into the
expression for $a_\tau$ gives that
\begin{equation*}
a_\tau=\O{|\Gamma_0|+\beta^2 A}
\end{equation*}
and hence
\begin{equation}
\ip{\xi_\tau}{\phi_{2a}}=\O{\beta^2 M_1(|\Gamma_0|+\beta^2 A) }.
\label{est:TimeDerxi}
\end{equation}

We now estimate the terms involving the linear operator ${\cal L}_{a
b c}$.  Write the operator ${\cal L}_{a b c}$ as
\begin{equation*}
{\cal L}_{a b c}={\cal L}_*+\frac{1}{4} a_\tau y^2-\frac{2 p
c}{p-1+b y^2},
\end{equation*}
where ${\cal L}_*$ is self-adjoint and satisfies ${\cal
L}_*\phi_{0a}=\frac{2 a}{p-1}\phi_{0a}$ and ${\cal L}_*\phi_{2a}=\frac{2 a p}{p-1}\phi_{2,a}$. Projecting
${\cal L}_{a b c}\xi$ onto the eigenvectors $\phi_{0a}$ and
$\phi_{2a}$ of ${\cal L}_*$ gives the equations
\begin{align*}
|\ip{{\cal L}_{a b
c}\xi}{\phi_{0a}}|&\lesssim | a_\tau| |\ip{\xi}{a y^2
e^{-\frac{a}{4} y^2}}|+|\ip{\xi}{\frac{by^{2}}{p-1+b
y^2} e^{-\frac{a}{4} y^2}}|=|
\ip{\xi}{\frac{by^{2}}{p-1+b y^2}
e^{-\frac{a}{4} y^2}}|,\\
|\ip{{\cal L}_{a b
c}\xi}{\phi_{2a}}|&=|a_\tau| |\ip{\xi}{a y^2(a
y^2-1)e^{-\frac{a}{4} y^2}}|+|\ip{\xi}{\frac{(a
y^2-1)by^{2}}{p-1+b y^2}e^{-\frac{a}{4} y^2}}|.
\end{align*}
Estimating with H\"{o}lder's inequality gives the estimates
\begin{align*}
|\ip{{\cal L}_{a b
c}\xi}{\phi_{0a}}|&\lesssim b\|\lra{y}^{-3}\xi e^{\frac{a}{4} y^2}\|_\infty\\
|\ip{{\cal L}_{a b c}\xi}{\phi_{2a}}|&\lesssim
(|a_\tau|+b)\|\lra{y}^{-3}\xi e^{\frac{a}{4} y^2}\|_\infty.
\end{align*}
In terms of the estimating functions $\beta$ and $M_1$, these
estimates, after using the above estimate of $a_\tau$ and
simplifying in $a$ and $c$, become
\begin{align}
\ip{{\cal L}_{a b
c}\xi}{\phi_{0a}}&\lesssim \beta^{3}M_1\label{Est:LinearOpPsi0}\\
\ip{{\cal L}_{a b c}\xi}{\phi_{2a}}&\lesssim \beta^2 M_1\lb
\beta+|\Gamma_0|+\beta^2 A \rb.\label{Est:LinearOpPsi2}
\end{align}

Lastly, we estimate the inner products involving the nonlinearity.
Due to \eqref{eqn:69a}, both $\ip{{\cal N}}{\phi_{0a}}$ and
$\ip{{\cal N}}{\phi_{2a}}$ are estimated by $\O{\|\lra{y}^{-3}
e^{\frac{a}{4}y^2}\xi\|_\infty^2+\|\lra{y}^{-3}
e^{\frac{a}{4}y^2}\xi\|_\infty^p}.$  Writing this in terms of
$\beta$ and $M_1$ and simplifying gives the estimate
\begin{equation}
|\ip{{\cal N}}{\phi_{ia}}|\lesssim \beta^4 M_1^2+\beta^{2p}
M_1^p.\label{est:NonPsi0Psi2}
\end{equation}

Estimates \eqref{est:ForcingPsi0}-\eqref{Est:LinearOpPsi2} and
\eqref{est:NonPsi0Psi2} imply that $\Gamma_0+\Gamma_1=R_1$ and
$\Gamma_1=R_2$, where $R_1$ and $R_2$ are of the order
\begin{equation*}
\O{\beta (|\Gamma_0|+|\Gamma_1|)+\beta^3+\beta^2 M_1\lb
\beta+|\Gamma_0|+\beta^2 A \rb+\beta^4 M_1^2+\beta^{2p} M_1^{p}}.
\end{equation*}
By the facts that $\beta(\tau)\le b_{0}\ll 1$ and $A, M_1\leq
\beta^{-\frac{\kappa}{2}} $, we obtain the estimates
\begin{equation}
|\Gamma_0|+|\Gamma_1|\lesssim \beta^3+\beta^3 M_1(1+ A)+\beta^4
M_1^2+\beta^{2p} M_1^{p} \label{GammaEst}
\end{equation}
for the times $0\le \tau\le T$.
\end{proof}
Equations (~\ref{eqn:Festimate}) and (~\ref{GammaEst}) yield the
following corollary.
\begin{cor}
\begin{equation}\label{eq:Festimate}
\begin{array}{lll}
\|\lra{y}^{-n} e^{\frac{a}{4} y^2}{\cal F}\|_\infty&\lesssim
&\beta^{k_{n}}(\tau)[1+ M_1(1+ A)+ M_1^2+ M_1^{p}]
\end{array}
\end{equation} with $n=0,3$ and $k_{0}:=\min\{1,2p-1\},$ $k_{3}:=\min\{5/2,2p\}.$
\end{cor}
\begin{remark}
Equation \eqref{eqn:FluctuationGauged} for the unknowns $a$, $b$,
$c$ and $\xi$ is invariant under the transformation
\begin{equation*}
(a(\tau), b(\tau), c(\tau),\xi(\tau))\mapsto (\mu^2 a(\mu\tau),\mu^2
b(\mu\tau),\mu^2 c(\mu\tau),\mu^\frac{2}{p-1}\xi(\mu y,\mu^2\tau)).
\end{equation*}
This symmetry is related to the symmetry \eqref{rescale} of
\eqref{NLH}. Consequently, Equations (~\ref{eqn:bParameter}) and
(~\ref{eqn:cParameter}) have the same symmetry.
\end{remark}
\begin{remark}
Dynamical equations (~\ref{eqn:bParameter}) and
(~\ref{eqn:cParameter}) have static solutions $(b,c,\xi)=(0,0,0)$
and $(b,c,\xi)=(0,a,0)$ with $a$ a constant (the latter implies
$a=\frac{1}{2}$).
\end{remark}

\section{Proof of Estimates
\label{Section:9} \eqref{EstB}-\eqref{EstA}}\label{SEC:EstB} Recall
that $a=2c-\frac{1}{2}$.
Assume $B(\tau)\le \beta^{-\frac{\kappa}{2}}(\tau)$ for
$\tau\in[0,T]$ which implies that $b\lesssim \beta,$
$\frac{1}{b}\lesssim\frac{1}{\beta}.$

We rewrite equation \eqref{eqn:bParameter} as $
b_\tau=-\frac{4p}{(p-1)^{2}} b^2+b\lb \frac{1}{2}-a-\frac{2 b}{p-1}
\rb+{\cal R}_b$.  By the definition of $A$, the second term on the
right hand side is bounded by $b \beta^{2} A\lesssim \beta^{3}A$.
Thus, using the bound for ${\cal R}_b$ given in Proposition
\ref{Prop:LyapunovSchmidtReduction}, we obtain \eqref{EstB}.

To prove \eqref{MajorE} we begin by dividing \eqref{EstB} by $b^2$
and using the inequality $\frac{1}{b}\lesssim \frac{1}{\beta}$ to
obtain the estimate
\begin{equation}
\left| -\p_\tau\frac{1}{b}+\frac{4p}{(p-1)^{2}} \right|\lesssim
\beta+\beta M_1(1+ A)+\beta^{2}M_{1}^{2}+\beta^{2p-2}M_{1}^{p}.
\label{est:InversebDE}
\end{equation}
Since $\beta$ is a solution to $-\p_\tau \beta^{-1}+4p(p-1)^{-2}=0$,
Equation \eqref{est:InversebDE} implies that
\begin{equation*}
|\p_\tau\lb \frac{1}{b}-\frac{1}{\beta}\rb|\lesssim \beta+\beta
M_1(1+ A)+\beta^{2}M_{1}^{2}+\beta^{2p-2}M_{1}^{p}.
\end{equation*}
Integrating this equation over $[0,\tau]$, multiplying the result by
$\beta^{-1-\kappa}$ and using that $\beta(0)=b(0)$, $b\lesssim
\beta$ gives the estimate
\begin{equation*}
\beta^{-1-\kappa}|\beta-b|\lesssim \beta^{1-\kappa} \int_0^\tau\lb
\beta+\beta M_1(1+ A)+\beta^{2}M_{1}^{2}+\beta^{2p-2}M_{1}^{p}\, \rb
ds,
\end{equation*}
where, recall, $\kappa:=\min\{\frac{1}{2}, \frac{p-1}{2}\}<1$.
Hence, by the definition of $\beta$ and $B$ and the facts that
$M_{1}$ and $A$ are increasing functions, \eqref{MajorE} follows.

Define the quantity $\Gamma:=\frac{1}{2}-a-\frac{2}{p-1} b$.
 Differentiating $\Gamma$ with respect to $\tau$ and substituting for
$b_\tau$ and $a_{\tau}=2c_\tau$ Equations \eqref{eqn:bParameter} and
\eqref{eqn:cParameter} we obtain
\begin{equation*}
\p_\tau\Gamma=-2 c(\Gamma+{\cal R}_c)-\frac{2}{p-1}\lb
-\frac{2(3p-1)}{(p-1)^2}b^2+2 b(c-a)+{\cal R}_b \rb.
\end{equation*}
Replacing $2b(c-a)$ by $b\Gamma+\frac{2}{p-1} b^2$ and rearranging
the resulting equation gives that
\begin{equation*}
\p_\tau \Gamma+\lsb a+\frac{1}{2}-\frac{2}{p-1} b \rsb\Gamma=\frac{8
p}{(p-1)^3} b^2-(a+\frac{1}{2}){\cal R}_c-\frac{2}{p-1}{\cal R}_b.
\end{equation*}
Let $\mu=\exp\lb \int_{0}^{\tau}\lb a+\frac{1}{2}-\frac{2}{p-1} b\rb
ds\rb$. Then the above equation implies that
\begin{equation*}
\p_\tau(\mu\Gamma)=\frac{8p}{(p-1)^3} \int_0^\tau \mu b^2\,
ds- \int_0^\tau (a+\frac{1}{2})\mu{\cal R}_c\, ds- \int_0^\tau
\frac{2}{p-1}\mu {\cal R}_b\, ds.
\end{equation*}
We now integrate the above equation over $[0,\tau]\subseteq [0,T]$
and use the inequality $b\lesssim \beta$ and the estimates of ${\cal
R}_b$ and ${\cal R}_c$ in Proposition
~\ref{Prop:LyapunovSchmidtReduction} to obtain
\begin{equation*}
|\Gamma|\lesssim \mu^{-1}\Gamma(0)+\mu^{-1} \int_0^\tau \mu \beta^2\,
ds+\mu^{-1} \int_0^\tau \mu\lb \beta^3+\beta^3 M_1(1+
A)+\beta^{4}M_{1}^{2}+\beta^{2p}M_{1}^{p} \rb\, ds.
\end{equation*}
For our purpose, it is sufficient to use the less sharp inequality
\begin{equation*}
|\Gamma|\lesssim\mu^{-1}\Gamma(0)+\lb
1+\beta(0)M_1(1+A)+\beta(0)M_{1}^{2}+\beta^{2p-2}(0)M_{1}^{p} \rb
\mu^{-1} \int_0^\tau \mu \beta^2\, ds.
\end{equation*} The assumption that $A(\tau),\ B(\tau)\leq
\beta^{-\frac{\kappa}{2}}(\tau)$ implies that
$a+\frac{1}{2}-\frac{2}{p-1}b=1-\frac{4 b}{p-1}+\O{\beta^2 A}\geq
\frac{1}{2}$ and therefore $\beta^{-2}\mu^{-1}\lesssim
\beta^{-2}(0)$ and $ \int_0^\tau \mu(s) \beta^2(s)\, ds\lesssim
\mu(\tau)\beta^2(\tau)$. The last two inequalities and the relation
$\displaystyle\max_{s\leq \tau}\beta^{-2}(s)|\Gamma(s)|=A(\tau)$
lead to \eqref{EstA}.
\section{Rescaling of Fluctuations on a Fixed Time Interval}
\label{Section:Rescaling} The coefficient in front of $y^2$ in the
operator $\mathcal{L}_{abc}$, \eqref{eqn:LinOp}, is time dependent,
complicating the estimation of the semigroup generated by this
operator. In this section we introduce the new time and space
variables in such a way that the coefficient at $y^2$ in the new
operator is constant (cf \cite{BP1, BuSu, Per}).

Let $T$ be given and let $t(\tau)$ be the inverse of the function
$\tau(t):= \int_0^t\lambda^2(s)\, ds$.
 We approximate the scaling parameter $\lambda(t)$ over the time
interval $[0,t(T)]$ by a new parameter $\lambda_1(t)$.  We choose
$\lambda_1(t)$ to satisfy for $t\leq t(T)$
\begin{equation*}
\p_t\lb \lambda_1^{-3}\p_t\lambda_1 \rb=0\ \mbox{with} \
\lambda_1(t(T))=\lambda(t(T))\ \mbox{and}\
\p_t\lambda_1(t(T))=\p_t\lambda(t(T)).
\end{equation*}
We define $\alpha:=\lambda_1^{-3}\p_t\lambda_1=a(T)$.  This is an
analog of the parameter $a$ and it is constant.  The last two
conditions imply that $\lambda_1$ is tangent to $\lambda$ at
$t=t(T)$.  Define the new time and space variables as
\begin{equation*}
z=\frac{\lambda_1}{\lambda} y\ \mbox{and}\ \sigma=\sigma(t(\tau))\
\text{with}\ \sigma(t):= \int_0^t \lambda_1^2(s)\, ds
\end{equation*}
where $\tau\leq T$, $\sigma\leq S:=\sigma(T)$ and $\lambda$
$\lambda_{1}$ are functions of $t(\tau).$ Now we introduce the new
function $\eta(z,\sigma)$ by the equality
\begin{equation}\label{NewFun}
\lambda_1^\frac{2}{p-1} e^{\frac{\alpha}{4}
z^2}\eta(z,\sigma)=\lambda^\frac{2}{p-1} e^{\frac{a}{4}
y^2}\xi(y,\tau).
\end{equation}

Denote by $t(\sigma)$ the inverse of the function $\sigma(t)$. In
the equation for $\eta(z,\sigma)$ derived below and in what follows
the symbols $\lambda$, $a$ and $b$ stand for $\lambda(t(\sigma)),$
$a(\tau(t(\sigma)))$ and $b(\tau(t(\sigma)))$, respectively.
 Substituting this change of variables into
\eqref{eqn:FluctuationGauged} gives the governing equation for
$\eta$:
\begin{equation}\label{eq:eta}
\p_\sigma\eta=-L_{\alpha}\eta+W(a,b,\alpha)\eta+F(a,b,\alpha)+N(\eta,a,
b,\alpha),
\end{equation}
where
\begin{align}
&L_{\alpha}:=L_0+V,\ L_0:=-\p_z^2+\frac{\alpha^2}{4}
z^2-\frac{5}{2}\alpha,\ V:=\frac{2 p\alpha}{p-1}-\frac{2 p
\alpha}{p-1+ \beta z^2}, \label{eqn:DefnL0V}\\
&W(a,b,\alpha):=\frac{\lambda^2}{\lambda_1^2}\frac{ p
(a+\frac{1}{2})}{p-1+b
\frac{\lambda^2}{\lambda_1^2} z^2}-\frac{2 p\alpha}{p-1+ \beta z^2},\nonumber\\
&F(a, b, \alpha):=\lb\frac{\lambda}{\lambda_1}\rb^\frac{2p}{p-1}
e^{-\frac{\alpha}{4} z^2}e^{\frac{a}{4}y^{2}}{\cal
F}(a,b,c)\nonumber
\end{align}
and
\begin{align*}
&N(\eta,a,b,\alpha):=\lb\frac{\lambda}{\lambda_1}\rb^\frac{2p}{p-1}
e^{-\frac{\alpha}{4}
z^2}e^{\frac{a}{4}\frac{\lambda^2}{\lambda_1^2}z^2}{\cal N}\lb
 \lb\frac{\lambda_1}{\lambda}\rb^\frac{2}{p-1} e^{\frac{\alpha}{4} z^2}
 e^{-\frac{a}{4}y^{2}}\eta,b,c\rb,
\end{align*}
where, recall, $c$ and $a$ are related as $2c=a+\frac{1}{2}$ and
$\beta$ is defined in \eqref{FunBTau}.

In the next statement we prove that the new parameter
$\lambda_{1}(t)$ is a good approximation of the old one,
$\lambda(t)$. We have
\begin{prop}\label{NewTrajectory}
If $A(\tau)\leq \beta^{-\frac{\kappa}{2}}(\tau)$ and $b(0)\ll 1$,
then
\begin{equation}\label{eq:appro}
|\frac{\lambda}{\lambda_{1}}(t(\tau))-1|\lesssim \beta(\tau)\leq
b(0).
\end{equation}
\end{prop}
\begin{proof}
Differentiating $\frac{\lambda}{\lambda_1}-1$ with respect to $\tau$
(recall that $\frac{dt}{d\tau}=\frac{1}{\lambda^2}$) gives the
expression
\begin{equation*}
\frac{d}{d\tau}\lb\frac{\lambda}{\lambda_1}-1\rb=\frac{\lambda}{\lambda_1}a-\frac{\lambda_1}{\lambda}\alpha
\end{equation*} or, after some manipulations
\begin{equation}\label{EstLambda}
\frac{d}{d\tau}[\frac{\lambda}{\lambda_{1}}-1]
=2a(\frac{\lambda}{\lambda_{1}}-1)+\Gamma
\end{equation} with
$$\Gamma:=a-\alpha-a\frac{\lambda_{1}}{\lambda}(\frac{\lambda}{\lambda_{1}}-1)^{2}+(a-\alpha)(\frac{\lambda_{1}}{\lambda}-1).$$
Observe that $\frac{\lambda}{\lambda_{1}}(t(\tau))-1=0$ when
$\tau=T.$ Thus Equations (~\ref{EstLambda}) can be rewritten as
\begin{equation}
\frac{\lambda}{\lambda_{1}}(t(\tau))-1=- \int_{\tau}^{T}e^{- \int^{\sigma}_{\tau}2a(\rho)d\rho}\Gamma(\sigma)d\sigma.
\label{eqn:89}
\end{equation}
By the definition of $A(\tau)$ and the definition $\alpha=a(T)$ we
have that, if $A(\tau)\leq \beta^{-\frac{\kappa}{2}}(\tau)$, then
\begin{equation}\label{CauchA}
|a(\tau)-\alpha|,\ |a(\tau)-\frac{1}{2}|\leq 2\beta(\tau)
\end{equation} on the time interval $\tau\in [0,T]$. Thus
\begin{equation}\label{Rem}
|\Gamma|\lesssim
\beta+(1+\frac{\lambda_{1}}{\lambda})(\frac{\lambda}{\lambda_{1}}-1)^{2}+\beta|\frac{\lambda}{\lambda_{1}}-1|.
\end{equation}
which together with \eqref{eqn:89} and \eqref{CauchA} implies
\eqref{eq:appro}.
\end{proof}
\section{Estimate on the Propagators}
\label{Section:PropEst} Let $\bar{P}^{\alpha}$ be the projection
onto the space spanned by the first three eigenvectors of $L_{0}$
and $P^{\alpha}:=1-\bar{P}^{\alpha}$. Denote by
$U^{(1)}_{\alpha}(\tau,\sigma)$ the propagator generated on ${\rm
Ran}\, P^\alpha$ by the operator $-P^{\alpha}L_{\alpha}P^{\alpha},$
where, recall, the definition of the operator $L_{\alpha}$ is given
in Equation (~\ref{eqn:DefnL0V}).
\begin{prop}\label{ProP}
For any function $g\in Ran P^{\alpha}$ and for $c_0:=\alpha -
\epsilon$ with some $\epsilon>0$ small we have
$$\|\langle z\rangle^{-3}e^{\frac{\alpha z^{2}}{4}}U^{(1)}_{\alpha}(\tau,\sigma)g\|_{\infty}\lesssim e^{-\cO(\tau-\sigma)}\|\langle z\rangle^{-3}e^{\frac{\alpha z^{2}}{4}}g\|_{\infty}.$$
\end{prop}
The proof of this proposition is given after Lemma ~\ref{LM:FK}.
Here we just observe that in the $L^{2}$ norm
$P^{\alpha}L_{\alpha}P^{\alpha}\geq
(-\partial_{z}^{2}+\frac{\alpha^{2}}{4}z^{2}-\frac{5}{2}\alpha)P^{\alpha}\geq
\frac{1}{2}\alpha P^{\alpha}$. However, this does not help in
proving the weighted $L^{\infty}$ bound above. We start with an
estimate for the propagator $U_{\alpha }(\tau,\sigma),$ generated by
the operator $-L_{\alpha}$. Recall the definition of the operator
$L_{0}$ in (~\ref{eqn:DefnL0V}) and define $U_0(x,y)$ as the
integral kernel of the operator $e^{-\frac{\alpha
z^{2}}{4}}e^{-rL_{0}}e^{\frac{\alpha z^{2}}{4}}$. We begin with
\begin{lemma}\label{kernelEst}
For $n=0,1,2,3,4$, any function $g$ and $r>0$ we have that
\begin{equation}
\|\langle z\rangle^{-n}e^{\frac{\alpha
z^{2}}{4}}e^{-L_{0}r}g\|_{\infty}\lesssim e^{2\alpha r}\|\langle
z\rangle^{-n}e^{\frac{\alpha z^{2}}{4}}g\|_{\infty} \label{est:99a}
\end{equation}
 or equivalently
\begin{equation}\label{eq:secondForm}
e^{\frac{\alpha x^{2}}{2}} \int \langle
x\rangle^{-n}U_{0}(x,y)e^{-\frac{\alpha}{2}y^{2}}\langle
y\rangle^{n}dy\lesssim e^{2\alpha r}.
\end{equation}
\end{lemma}
\begin{proof} We only prove the case $n=2.$ The cases $n=0,4$ are
similar.  The cases $n=1,3$ follows from $n=0,2,4$ by an
interpolation result. Note that the first four eigenvectors of
$L_{0}$ are $e^{-\frac{\alpha x^{2}}{4}},\ xe^{-\frac{\alpha
x^{2}}{4}},\ (\alpha x^{2}-1)e^{-\frac{\alpha x^{2}}{4}}$ and
$(\alpha x^{3}-3x)e^{-\frac{\alpha x^{2}}{4}}$ with the eigenvalues
$-2\alpha,\ -\alpha,\ 0$ and $\alpha.$ Thus for the case $n=2$,
using that the integral kernel of $e^{-r L_0}$ is positive and
therefore $\|e^{-r L_0} g\|_\infty\le \|f^{-1} g\|_\infty\|e^{-r
L_0} f\|_\infty$ for any $f>0$ and using that $e^{-r L_0}
e^{-\frac{\alpha}{4} z^2}=e^{2\alpha r}e^{-\frac{\alpha}{4}z^2}$ and
$e^{-r L_0}(\alpha z^2-1)e^{-\frac{\alpha}{4} z^2}=(\alpha z^2-1)e^{-\frac{\alpha}{4} z^2}$, we find that
$$
\begin{array}{lll}
\|\langle z\rangle^{-2}e^{\frac{\alpha
z^{2}}{4}}e^{-rL_{0}}g\|_{\infty} &\leq &\|\langle z\rangle^{-2}
e^{\frac{\alpha z^{2}}{4}}e^{-rL_{0}}e^{-\frac{\alpha
z^{2}}{4}}(z^{2}+1)\|_{\infty}\|\langle z\rangle^{-2}e^{\frac{\alpha
z^{2}}{4}}g\|_{\infty}\\
&=& \|\langle z\rangle^{-2}[e^{2\alpha
r}\frac{1}{\alpha}+(z^{2}-\frac{1}{\alpha})]\|_{\infty}\|\langle
z\rangle^{-2}e^{\frac{\alpha
z^{2}}{4}}g\|_{\infty}\\
&\leq&2(\frac{1}{\alpha}+1)e^{2\alpha r}\|\langle
z\rangle^{-2}e^{\frac{\alpha z^{2}}{4}}g\|_{\infty}.
\end{array}
$$
This implies \eqref{est:99a}.  To prove \eqref{eq:secondForm} we
note that $U_0(x,y)$ is, by definition, the integral kernel of the
operator $e^{-\frac{\alpha}{4} z^2}e^{-r L_0} e^{\frac{\alpha}{4}
z^2}$.  Thus, taking $g(x)=\lra{x}^n e^{-\frac{\alpha}{4} x^2}$ in
\eqref{est:99a} yields \eqref{eq:secondForm}.
\end{proof}
A version of the following lemma is proved in \cite{BrKu}.
\begin{lemma}\label{propagator}
For any function $g$ and positive constants $\sigma$ and $r$ we have
$$\|\langle z\rangle^{-3}e^{\frac{\alpha z^{2}}{4}}U_{\alpha}(\sigma+r,\sigma)P^{\alpha}g\|_{\infty}
\lesssim [e^{2\alpha r}r(1+r){\beta^{1/2}(\sigma)}+e^{-\alpha
r}]\|\langle z\rangle^{-3}e^{\frac{\alpha z^{2}}{4}}g\|_{\infty}.$$
\end{lemma}
\begin{proof}
The spatial variables in this proof will be denoted by $x$, $y$ and
$z$.  Recall the definitions of the operators $L_0$ and $V$ in
\eqref{eqn:DefnL0V}. Denote the integral kernel of $e^{-\frac{\alpha
x^{2}}{4}}U_{\alpha}(\sigma+r,\sigma)e^{\frac{\alpha x^{2}}{4}}$ by
$U(x,y)$. By Theorem ~\ref{THM:trotter}, given in Appendix B below,
we have the representation
\begin{equation}\label{eq:FK1}
U(x,y)=U_{0}(x,y)\langle e^{V}\rangle(x,y),
\end{equation}
where, recall that $U_0(x,y)$ is the integral kernel of the operator
$e^{-\frac{\alpha z^{2}}{4}}e^{-rL_{0}}e^{\frac{\alpha z^{2}}{4}}$
and
\begin{equation}\label{FK2}
\langle e^{V}\rangle(x,y)= \int_{\sigma}^{\sigma+r}
e^{  \int_{\sigma}^{\sigma+r}-V(\sigma+s,\omega(s)+\omega_{0}(s))ds}
d\mu(\omega).
\end{equation}
Here $\omega_{0}(s)$ is defined in Theorem ~\ref{THM:trotter} of
Appendix \ref{Sec:Trotter} and $d\mu(\omega)$ is a harmonic
oscillator (Ornstein-Uhlenbeck) probability measure on the
continuous paths $\omega: [\sigma,\sigma+r]\rightarrow \mathbb{R}$
with the boundary condition $\omega(\sigma)=\omega(\sigma+r)=0.$ By
a standard formula (see \cite{Simon, GlJa}) we have
\begin{equation*}
U_{0}(x,y)=4\pi (1-e^{-2\alpha r})^{-1/2}\sqrt{\alpha}e^{2\alpha
r}e^{-\alpha\frac{(x-e^{-\alpha r}y)^{2}}{2(1-e^{-2\alpha r})}}.
\label{eqn:96a}
\end{equation*}

Define a new function $f:=e^{-\frac{\alpha y^{2}}{4}}P^{\alpha}g$.
The definitions above imply
\begin{equation}\label{FK1}
U_{\alpha}(\sigma+r,\sigma)P^{\alpha}g= \int e^{\frac{\alpha
x^{2}}{4}}U_{0}(x,y)\langle e^{V}\rangle(x,y)f(y)dy.
\end{equation}
Integrate by parts on the right hand side of (~\ref{FK1}) to obtain
\begin{equation}\label{Estimateta}
\begin{array}{lll}
U_{\alpha}(\sigma+r,\sigma)P^{\alpha}g&=&\displaystyle\sum_{k=0}^{2}e^{\frac{\alpha
x^{2}}{4}} \int
\partial_{y}^{k}U_{0}(x,y)\partial_{y}\langle
e^{V}\rangle(x,y)f^{(-k-1)}(y)dy\\
& &+e^{\frac{\alpha x^{2}}{4}} \int \partial_{y}^{3}U_{0}(x,y)\langle
e^{V}\rangle(x,y)f^{(-3)}(y)dy
\end{array}
\end{equation} where
$f^{(-m-1)}(x):= \int_{-\infty}^{x}f^{(-m)}(y)dy$ and $f^{(-0)}:=f.$
Now we estimate every term on the right hand side of Equation
(~\ref{Estimateta}).
\begin{enumerate}
\item[(A)] By the facts that $f=e^{-\frac{\alpha
y^{2}}{4}}P^{\alpha}g$ and $P^{\alpha}g\perp y^{n}e^{-\frac{\alpha
y^{2}}{4}},\ n=0,1,2,3,$ we have that $f\perp 1,\ y,\ y^{2},\
y^{3}.$ Therefore by integration by parts we have
$$f^{(-m)}(y)= \int_{-\infty}^{y}f^{(-m+1)}(x)dx=- \int_{y}^{\infty}f^{(-m+1)}(x)dx,\ m=1,2,3.$$
Moreover by the definition of $f^{(-m)}$ and the equation above we
have
$$|f^{(-m)}(y)|\lesssim\langle
y\rangle^{3-m}e^{-\frac{\alpha}{2}y^{2}}\|\langle
y\rangle^{-3}e^{\frac{\alpha}{4}y^{2}}P^{\alpha}g\|_{\infty}.$$
\item[(B)] Using the explicit formula for $U_0(x,y)$ given above we find $$|\partial^{(k)}_{y}U_{0}(x,y)|\lesssim \frac{e^{-\alpha
kr}}{(1-e^{-2\alpha r})^{k}}(|x|+|y|+1)^{k}U_{0}(x,y)$$.
\item[(C)] By an estimate from Appendix ~\ref{Sec:Trotter} (see also \cite{BrKu})  we have that
\begin{equation}\label{eqn:93a}
|\partial_{y}\langle e^{V}\rangle(x,y)|\leq {\beta^{1/2}}r.
\end{equation}
\end{enumerate}

Collecting the estimates (A)-(C) above and using Equation
(~\ref{Estimateta}), we have the following result
$$
\begin{array}{lll}
& &\langle x\rangle^{-3}e^{\frac{\alpha x^{2}}{4}}|U_{\alpha}(\sigma+r,\sigma)P^{\alpha}g(x)|\\
&\lesssim & \frac{\sqrt{\beta}r(1+r)}{(1-e^{-2\alpha r})^{3}}\langle
x\rangle^{-3}e^{\frac{\alpha
x^{2}}{2}}\displaystyle\sum_{k=0}^{2} \int
(|x|+|y|+1)^{k+1}U_{0}(x,y)|f^{(-k-1)}(y)|dy\\
& &+\frac{1}{(1-e^{-2\alpha r})^{3}}\langle
x\rangle^{-3}e^{\frac{\alpha x^{2}}{2}} \int (|x|+|y|+1)^{3}e^{-3\alpha
r}U_{0}(x,y)|f^{(-3)}(y)|dy\\
&\lesssim &\frac{\beta^{1/2}r(1+r)+e^{-3\alpha r}}{(1-e^{-2\alpha
r})^{3}}\sum_{n=0}^{3}e^{\frac{\alpha x^{2}}{2}} \int \langle
x\rangle^{-n}U_{0}(x,y)e^{-\frac{\alpha}{2}y^{2}}\langle
y\rangle^{n}dy \|\langle
y\rangle^{-3}e^{\frac{\alpha}{4}y^{2}}P^{\alpha}g\|_{\infty}.
\end{array}
$$
This together with the estimate \eqref{eq:secondForm} of Lemma
\ref{kernelEst} gives the estimate of Lemma \ref{propagator}.
\end{proof}
We will also need
\begin{lemma}\label{LM:FK}
\begin{equation}\label{Mehler}
\|\langle z\rangle^{-n}e^{\frac{\alpha
z^{2}}{4}}U_{\alpha}(\tau,\sigma)g\|_{\infty}\leq
e^{2\alpha(\tau-\sigma)}\|\langle z\rangle^{-n}e^{\frac{\alpha
z^{2}}{4}}g\|_{\infty}
\end{equation} with $n=0\ \text{or}\ 3.$
\end{lemma}
\begin{proof}
By Equations (~\ref{FK2}) and (~\ref{FK1}) we have that
$|U_{\alpha}(\tau,\sigma)|(x,y)\leq e^{-L_{0}(\tau-\sigma)}(x,y).$
Thus we have
\begin{equation}\label{FeyKac}
\begin{array}{lll}
\|\langle z\rangle^{-n}e^{\frac{\alpha
z^{2}}{4}}U_{\alpha}(\tau,\sigma)g\|_{\infty}& \le & \|\langle
z\rangle^{-n}e^{\frac{\alpha
z^{2}}{4}}e^{-L_{0}(\tau-\sigma)}|g|\|_{\infty}.
\end{array}
\end{equation}  Now we use Lemma ~\ref{kernelEst} to
estimate the right hand side to complete the proof.
\end{proof}

{\textbf{Proof of Proposition ~\ref{ProP}}}.  Recall that
$\bar{P}_\alpha$ is the projection on the span of the three first
eigenfunctions of the operator $L_0$ and
$P^{\alpha}:=1-\bar{P}^{\alpha}$. We write
\begin{equation}\label{TranForm}
L_{\alpha}=P^\alpha L_{\alpha} P^\alpha+E_1+\bar{P}^\alpha L_{\alpha
}\bar{P}^\alpha,
\end{equation}
where the operator $E_{1}$ is defined as $E_{1}:=\bar{P}^{\alpha}
L_{\alpha }P^\alpha+P^\alpha L_{\alpha }\bar{P}^\alpha.$ Using that
$\bar{P}^\alpha P^\alpha=0$, we transform $E_1$ to
$$
\begin{array}{lll}
E_{1}&=&-\bar{P}^{\alpha}\frac{2p(p-1)\alpha\beta z^{2}}{p-1+\beta
z^{2}}P^{\alpha}-P^{\alpha}\frac{2p(p-1)\alpha\beta z^{2}}{p-1+\beta
z^{2}}\bar{P}^{\alpha}.
\end{array}
$$  This implies
\begin{equation}\label{EstD2}
\|\langle z\rangle^{-3}e^{\frac{\alpha
z^{2}}{4}}E_{1}\eta(\sigma)\|_{\infty}\lesssim
\beta(\tau(\sigma))\|\langle z\rangle^{-3}e^{\frac{\alpha
z^{2}}{4}}\eta(\sigma)\|_{\infty}.
\end{equation}
We use Duhamel's principle to rewrite the propagator
$U^{(1)}_{\alpha}(\sigma_{1},\sigma_{2})$ on ${\rm Ran}\, P^\alpha$
as
\begin{equation}
U_{\alpha}^{(1)}(\sigma_{1},\sigma_{2})P^\alpha=
U_{\alpha}(\sigma_{1},\sigma_{2})P^\alpha- \int_{\sigma_{2}}^{\sigma_{1}}
U_{\alpha}(\sigma_{1},s)E_{1}U_{\alpha}^{(1)}(s,\sigma_{2}) P^\alpha
ds. \label{eqn:106}
\end{equation}
  Let $r=\sigma_{1}-\sigma_{2}$, $g\in {\rm Ran} P^{\alpha}$ and
$\eta(\sigma_{1}):=U^{(1)}_{\alpha}(\sigma_{1},\sigma_{2})g$.  We
estimate the two terms on the right hand side of \eqref{eqn:106}. We
claim that if $e^{ \alpha r}\leq \beta^{-1/32}(\tau(\sigma_{2}))$
then we have
\begin{equation}\label{firstIter}
\|\langle z\rangle^{-3}e^{\frac{\alpha
z^{2}}{4}}\eta(\sigma_{1})\|_{\infty}\lesssim e^{- \alpha
r}\|\langle z\rangle^{-3}e^{\frac{\alpha
z^{2}}{4}}\eta(\sigma_{2})\|_{\infty}.
\end{equation} To prove the claim we compute each terms on the right
hand side of (~\ref{firstIter}).
\begin{itemize}
\item[(A)] Notice that $P^{\alpha}\eta(s)=\eta(s)$.  We use Lemma
~\ref{propagator} to obtain, for $e^{\alpha r}\leq
\beta^{-1/32}(\tau(\sigma_{2}))$,
\begin{equation}\label{SecTerm} \|\langle
z\rangle^{-3}e^{\frac{\alpha
z^{2}}{4}}U_{\alpha}(\sigma_{1},\sigma_{2})g\|_{\infty}\lesssim
e^{-\alpha r}\|\langle z\rangle^{-3}e^{\frac{\alpha
z^{2}}{4}}g\|_{\infty}.\end{equation}

\item[(B)]
By Lemma ~\ref{LM:FK} and (~\ref{EstD2}) we obtain
\begin{equation*}
\begin{array}{lll}
\|\langle z\rangle^{-3}e^{\frac{\alpha
z^{2}}{4}} \int_{\sigma_{2}}^{\sigma_{1}}U_{\alpha}(\sigma_{1},s)E_{1}\eta(s)ds\|_{\infty}
\lesssim
 \int_{\sigma_{2}}^{\sigma_{1}}e^{2\alpha(\sigma_{1}-s)}\beta(\tau(s))\|\langle
z\rangle^{-3}e^{\frac{\alpha z^{2}}{4}}\eta(s)\|ds.
\end{array}
\end{equation*}
Using the condition $e^{\alpha r}\leq \beta^{-1/32}(\sigma_{2})$ and
the relation $\beta(\tau(s))\leq \beta(\tau(\sigma_{2}))$ for $s\geq
\sigma_{2}$ again, we find
\begin{equation}\label{ThirTerm}
\|\langle
z\rangle^{-3}e^{\frac{\alpha z^{2}}{4}} \int_{\sigma_{2}}^{\sigma_{1}}U_{\alpha}(\sigma_{1},s)E_{1}\eta(s)ds\|_{\infty}\\
\lesssim
 \int_{\sigma_{2}}^{\sigma_{1}}e^{-\alpha(\sigma_{1}-s)}\beta^{1/2}(\tau(s))\|\langle
z\rangle^{-3}e^{\frac{\alpha z^{2}}{4}}\eta(s)\|ds.
\end{equation}
\end{itemize}
Equations (~\ref{eqn:106}), (~\ref{SecTerm}) and (~\ref{ThirTerm})
imply that if $e^{\alpha r}\leq \beta^{-1/32}(\tau(\sigma_{2}))$
then (remember that $\eta(\sigma_{2})=g$)
\begin{equation}\label{FinalTerm}
\begin{array}{lll}
\|\langle z\rangle^{-3}e^{\frac{\alpha
z^{2}}{4}}\eta(\sigma_{1})\|_{\infty}&\lesssim & e^{-\alpha
r}\|\langle z\rangle^{-3}e^{\frac{\alpha
z^{2}}{4}}\eta(\sigma_{2})\|_{\infty}
+ \int_{\sigma_{2}}^{\tau}e^{-\alpha(\tau-s)}\beta^{1/2}(\tau(s))\|\langle
z\rangle^{-3}e^{\frac{\alpha z^{2}}{4}}\eta(s)\|ds.
\end{array}
\end{equation}

Next, we define a function $K(r)$ as
\begin{equation}\label{defineKz}
K(r):=\max_{0\leq k\leq r}e^{\alpha k} \|\langle
z\rangle^{-3}e^{\frac{\alpha z^{2}}{4}}\eta(\sigma_{2}+k)\|.
\end{equation}
Then (~\ref{FinalTerm}) implies that $$K(\sigma_{1})\lesssim
\|\langle z\rangle^{-3}e^{\frac{\alpha
z^{2}}{4}}\eta(\sigma_{2})\|_{\infty}+ \int_{\sigma_{2}}^{\sigma_{1}}e^{-\alpha(\sigma_{1}-s)}e^{-\alpha(s-\sigma_{2})}\beta^{1/2}(\tau(s))ds
K(\sigma_{2}).$$ We observe that
$$ \int_{\sigma_{2}}^{\sigma_{1}}e^{-\alpha(\sigma_{1}-s)}e^{-\alpha(s-\sigma_{2})}\beta^{1/2}(\tau(s))ds\leq
1/2$$ if $\beta(0)$ and, therefore, $\beta(\tau(s))=
\frac{1}{\frac{1}{\beta(0)}+\frac{4p}{(p-1)^{2}} \tau(s)}$ are
small. Thus we have
$$K(\sigma_{1})\lesssim \|\langle
z\rangle^{-3}e^{\frac{\alpha
z^{2}}{4}}\eta(\sigma_{2})\|_{\infty},$$ which together with
Equation (~\ref{defineKz}) implies (~\ref{firstIter}). Iterating
(~\ref{firstIter}) completes the proof of the proposition.
\begin{flushright}
$\square$
\end{flushright}

\section{Estimate of $M_{1}(\tau)$ (Equation \eqref{M1})}\label{SEC:EstM1} In this subsection we derive an
estimate for $M_{1}(T)$ given in Equation (~\ref{M1}).  Given any
time $\tau^{'}$, choose $T=\tau^{'}$ and pass from the unknown
$\xi(y,\tau)$, $\tau\leq T,$ to the new unknown $\eta(z,\sigma),$
$\sigma\leq S,$ given in (~\ref{NewFun}). Now we estimate the latter
function. To this end we use Equation (~\ref{eq:eta}). Observe that
the function $\eta$ is not orthogonal to the first three
eigenvectors of the operator $L_{0}$ defined in
(~\ref{eqn:DefnL0V}). Thus we apply the projection $P^\alpha$ to
Equation (~\ref{eq:eta}) to get
\begin{equation}\label{EQ:eta2}
\frac{d}{d\sigma}P^{\alpha}\eta=-P^{\alpha}L_{\alpha}P^{\alpha}\eta+P^{\alpha}\sum_{n=1}^{4}D_{n},
\end{equation} where
we used the fact that $P^\alpha$ are $\tau$-independent and the
functions $D_{n}\equiv D_{n}(\sigma), \ n=1,2,3,4,$ are defined as
$$D_{1}:=-P^{\alpha}V\eta+P^{\alpha}VP^{\alpha}\eta,\ \ \ D_{2}:=W(a,b,\alpha)\eta,$$
$$D_{3}:=F(a,b,\alpha),\ \  D_{4}:=N(\eta,a,b,\alpha),$$ recall the definitions of the functions $V$, $W$, $F$
and $N$ after \eqref{eqn:DefnL0V}.
\begin{lemma}\label{LM:EstDs} If $A(\tau),\ B(\tau)\leq
\beta^{-\frac{\kappa}{2}}(\tau)$ for $\tau\leq T$ and $b_0\ll 1$,
then we have
\begin{equation}\label{eq:estD1}
\|\langle z\rangle^{-3}e^{\frac{\alpha
z^{2}}{4}}D_{1}(\sigma)\|_{\infty}\lesssim
\beta^{5/2}(\tau(\sigma))M_{1}(T),
\end{equation}
\begin{equation}\label{eq:Remainder3}
\|\langle z\rangle^{-3}e^{\frac{\alpha
z^{2}}{4}}D_{2}(\sigma)\|_{\infty}\lesssim{\beta^{2+\frac{\kappa}{2}}(\tau
(\sigma))}M_{1}(T),
\end{equation}
\begin{equation}\label{eq:Festimate2}
\|\lra{z}^{-3} e^{\frac{\alpha}{4} z^2} D_{3}(\sigma)\|_\infty
\lesssim\beta^{\min\{5/2,2p\}}(\tau(\sigma))[1+ M_1(T)(1+ A(T))+
M_1^2(T)+ M_1^{p}(T)],
\end{equation}
\begin{equation}\label{nonlinearity3}
\|\langle z\rangle^{-3}e^{\frac{\alpha
}{4}z^{2}}D_{4}\|_{\infty}\lesssim
\beta^{2}(\tau(\sigma))M_{1}(T)[\beta^{1/2}(\tau(\sigma))M_{1}(T)+M_{2}(T)+\beta^{\frac{p-1}{2}}(\tau(\sigma))M_{1}^{p-1}(T)+M_{2}^{p-1}(T)].
\end{equation}
\end{lemma}
\begin{proof}
In what follows we use the following estimates, implied by
(~\ref{eq:appro}),
\begin{equation}\label{eq:compare2}
\frac{\lambda_{1}}{\lambda}(t(\tau))-1=O(\beta(\tau)), \
\text{thus}\ \frac{\lambda_{1}}{\lambda}(t(\tau)),\
\frac{\lambda}{\lambda_{1}}(t(\tau))\leq 2,\ \langle
z\rangle^{-3}\lesssim \langle y\rangle^{-3}
\end{equation} where, recall that $z:=\frac{\lambda_{1}}{\lambda}y.$
We start with proving the following two estimates which will be used
frequently below
\begin{equation}\label{Compare0}
\| e^{\frac{\alpha z^{2}}{4}}\eta(\sigma)\|_{\infty}\lesssim
\beta^{1/2}(\tau(\sigma))M_{1}(T)+M_{2}(T),
\end{equation}
\begin{equation}\label{Compare3}
\|\langle z\rangle^{-3} e^{\frac{\alpha
z^{2}}{4}}\eta(\sigma)\|_{\infty}\lesssim
\beta^{2}(\tau(\sigma))M_{1}(T).
\end{equation}
%
%
Recall the definition of $\chi_{\geq D}$ from (~\ref{cutoff}).
Writing  $1= 1 - \chi_{\geq D} + \chi_{\geq D}$ and using the
inequality $1-\chi_{\geq D}\lesssim \beta^{-3/2}(\tau)\langle
y\rangle^{-3}$, the relation between $\xi$ and $\eta$, see
(~\ref{NewFun}), and Estimate (~\ref{eq:compare2}) we find
\begin{equation}\label{eq:keyEst}
\begin{array}{lll}
\|e^{\frac{\alpha z^{2}}{4}}\eta(\sigma)\|_{\infty}\lesssim
\|e^{\frac{a(\tau(\sigma))y^{2}}{4}}\xi(\tau(\sigma))\|_{\infty}\lesssim
 \beta^{-3/2}(\tau(\sigma))\|\langle
y\rangle^{-3}e^{\frac{a(\tau(\sigma))y^{2}}{4}}\xi(\tau(\sigma))\|_{\infty}\\+
\|e^{\frac{a}{4}y^{2}} \chi_{\geq D}\xi(\tau)\|_{\infty}\ \leq
\beta^{1/2}(\tau(\sigma))M_{1}(T)+M_{2}(T)
\end{array}
\end{equation}
which is (~\ref{Compare0}). Similarly recall that
$z=\frac{\lambda_{1}}{\lambda}y$ which together with (~\ref{NewFun})
and (~\ref{eq:compare2}) yields
$$\|\langle z\rangle^{-3}e^{\frac{\alpha z^{2}}{4}}\eta(\sigma)\|_{\infty}\lesssim
\|\langle
y\rangle^{-3}e^{\frac{a(\tau(\sigma))y^{2}}{4}}\xi(\tau(\sigma))\|_{\infty}\lesssim
\beta^{2}(\tau(\sigma))M_{1}(T).$$ Thus we have (~\ref{Compare3}).

Now we proceed directly to proving the lemma. First we rewrite
$D_{1}$ as
$$D_{1}(\sigma)=-P^{\alpha}\frac{2p\alpha}{(p-1)(p-1+\beta(\tau(\sigma))z^{2})}\beta(\tau(\sigma))z^{2}(1-P^{\alpha})\eta(\sigma).$$
Now, using that $\langle z\rangle^{-1}\frac{b z^{2}}{1+b
z^{2}}\lesssim b^{1/2}$ and that $b\lesssim \beta$, we obtain
$$
\begin{array}{lll}
\|\langle z\rangle^{-3}e^{\frac{\alpha
z^{2}}{4}}D_{1}(\sigma)\|_{\infty}&\lesssim& \beta^{1/2}(\tau)|\|
\langle z\rangle^{-2}e^{\frac{\alpha
z^{2}}{4}}(1-P^{\alpha})\eta(\sigma)\|_{\infty}.
\end{array}
$$ Next, due to the explicit
form of $\bar{P}^{\alpha}:=1-P^{\alpha}$, i.e.
$\bar{P}^{\alpha}=\displaystyle\sum_{m=0}^{2}| \phi_{m,\alpha}
\rangle \langle \phi_{m,\alpha} |,$
where $\phi_{m,\alpha}$ are the normalized eigenfunctions of the
operator $L_0:=-\p_z^2+\frac{\alpha^2}{4} z^2-\frac{5}{2}\alpha,$
and decay properties of these eigenfunctions, see
(~\ref{eq:eigenvectors}) of Appendix \ref{Appendix:SpectrumLinear}
below,
we have for any function $g$
\begin{equation}\label{eq:estPro} \|\langle
z\rangle^{-2}e^{\frac{\alpha z^{2}}{4}}\bar{P}^{\alpha}
g\|_{\infty}\lesssim \|\langle z\rangle^{-3}e^{\frac{\alpha
z^{2}}{4}}g\|_{\infty}.
\end{equation} Collecting
the estimates above and using (~\ref{Compare3}), we arrive at
$$\|\langle z\rangle^{-3}e^{\frac{\alpha
z^{2}}{4}}D_{1}(\sigma)\|_{\infty}\lesssim
\beta^{1/2}(\tau(\sigma))\|\langle z\rangle^{-3}e^{\frac{\alpha
z^{2}}{4}}\eta(\sigma)\|_{\infty} \lesssim
\beta^{5/2}(\tau(\sigma))M_{1}(T).$$

To prove (~\ref{eq:Remainder3}) we recall the definition of $D_{2}$
and rewrite it as
$$
\begin{array}{lll}
D_{2}&=&\{[\frac{\lambda^{2}}{\lambda_{1}^{2}}-1]p\frac{2a+1}{p-1+by^{2}}+p\frac{2(a-\alpha)}{p-1+by^{2}}+\frac{b(2\alpha+1)(\frac{\lambda_{1}^{2}}{\lambda^{2}}-1)y^{2}}{(p-1+b
z^{2})(p-1+b y^{2})}+p\frac{1-2\alpha}{p-1+b z^{2}}+\frac{2p\alpha
\beta z^{2}}{(p-1+bz^{2})(p-1+\beta
z^{2})}\frac{\beta-b}{\beta}\}\eta.
\end{array}
$$
Then Equations (~\ref{eq:appro}), (~\ref{CauchA}) and the definition
of $B$ in (~\ref{majorants}) imply
$$\|\langle z\rangle^{-3}e^{\frac{\alpha}{4} z^2}D_{2}(\sigma)\|_{\infty}\leq \beta^{\frac{\kappa}{2}}(\tau(\sigma))\|\langle y\rangle^{-3}e^{\frac{\alpha z^{2}}{4}}\eta(\sigma)\|_{\infty}.$$
Using (~\ref{Compare3}) we obtain (~\ref{eq:Remainder3}) (recall
$\kappa:=\min\{\frac{1}{2},\frac{p-1}{2}\}$).

Now we prove (~\ref{eq:Festimate2}). By (~\ref{eq:compare2}) and the
relation between $D_{3}$, $F$ and $\mathcal{F}$ we have
$$\|\langle z\rangle^{-3}e^{\frac{\alpha}{4} z^2} D_{3}(\sigma)\|_\infty \lesssim \|\langle y\rangle^{-3}e^{\frac{a(\tau(\sigma))}{4} y^2}
\mathcal{F}(a,b,c)(\tau(\sigma))\|_\infty$$ which together with
(~\ref{eq:Festimate}) implies (~\ref{eq:Festimate2}).

Lastly we prove (~\ref{nonlinearity3}). By the relation between
$D_{4}$, $N$ and $\mathcal{N}$ and the estimate in (~\ref{eqn:69a})
we have
$$
\begin{array}{lll}
\|\langle z\rangle^{-3}e^{\frac{\alpha
}{4}z^{2}}D_{4}(\sigma)\|_{\infty}&\lesssim &\|\langle
y\rangle^{-3}e^{\frac{a(\tau(\sigma))}{4}y^{2}}\mathcal{N}(\xi(\tau(\sigma)),b(\tau(\sigma)),c(\tau(\sigma)))\|_{\infty}\\
&\lesssim &\|\langle y\rangle^{-3}e^{\frac{a
y^{2}}{4}}\xi(\tau(\sigma))\|_{\infty}[\|e^{\frac{a
y^{2}}{4}}\xi(\tau(\sigma))\|_{\infty}+\|e^{\frac{a
y^{2}}{4}}\xi(\tau(\sigma))\|_{\infty}^{p-1}].
\end{array}
$$ Using (~\ref{eq:keyEst}) and the definition of $M_{1}$ we complete
the proof.
\end{proof}

Below we will need the following lemma. Recall that
$S:=\sigma(t(T))$.
\begin{lemma}\label{Bridge}  If $A(\tau)\leq
\beta^{-\frac{\kappa}{2}}(\tau)$, then for any $c_{1},c_{2}>0$ there
exists a constant $c(c_{1},c_{2})$ such that
\begin{equation}\label{INT}
 \int_{0}^{S}e^{-c_{1}(S-\sigma)}\beta^{c_{2}}(\tau(t(\sigma)))d\sigma\leq
c(c_{1},c_{2})\beta^{c_{2}}(T).
\end{equation}
\end{lemma}
\begin{proof}
We use the shorthand $\tau(\sigma)\equiv \tau(t(\sigma)),$ where,
recall $t(\sigma)$ is the inverse of
$\sigma(t)= \int_{0}^{t}\lambda_{1}^{2}(k)dk$ and
$\tau(t)= \int_{0}^{t}\lambda^{2}(k)dk.$ By Proposition
~\ref{NewTrajectory} we have that $\frac{1}{2}\leq
\frac{\lambda}{\lambda_{1}}\leq 2$ provided that $A(\tau)\leq
\beta^{-\frac{\kappa}{2}}(\tau)$. Hence
\begin{equation}\label{TauS1}
\frac{1}{4}\sigma\leq \tau(\sigma)\leq 4\sigma
\end{equation} which
implies
$\frac{1}{\frac{1}{b(0)}+\frac{4p}{(p-1)^{2}}\tau(\sigma)}\lesssim
\frac{1}{\frac{1}{b(0)}+\sigma}$.  By a direct computation we have
\begin{equation}\label{InI2}
 \int_{0}^{S}e^{-c_{1}(S-\sigma)}\beta^{c_{2}}(\tau(\sigma))d\sigma\leq
c(c_{1},c_{2})\frac{1}{(\frac{1}{b(0)}+\frac{4p}{p-1}S)^{c_{2}}}.
\end{equation} Using
(~\ref{TauS1}) again we obtain $4S\geq \tau(S)=T\geq \frac{1}{4}S$
which together with (~\ref{InI2}) implies \eqref{INT}.
\end{proof}

Recall that $U_{\alpha}^{(1)}(t,s)$ is the propagator generated by
the operator $-P^{\alpha} L_{\alpha}P^{\alpha}$. To estimate the
function $P^{\alpha}\eta$ we rewrite Equation (~\ref{EQ:eta2}) as
$$
P^{\alpha}\eta(S)=U_{\alpha}^{(1)}(S,0)P^{\alpha}\eta(0)+\displaystyle\sum_{n=1}^{4} \int_{0}^{S}U^{(1)}_{\alpha}(S,\sigma)P^{\alpha}D_{n}(\sigma)d\sigma
$$ which implies
\begin{equation}\label{eq:estEsta}
\|\langle z\rangle^{-3}e^{\frac{\alpha
}{4}z^{2}}P^{\alpha}\eta(S)\|_{\infty}\leq K_{1}+K_{2}
\end{equation} with $$K_{1}:=\|\langle
z\rangle^{-3}e^{\frac{\alpha
z^{2}}{4}}U_{\alpha}^{(1)}(S,0)P^{\alpha}\eta(0)\|_{\infty};$$
$$K_{2}:=\|\langle
z\rangle^{-3}e^{\frac{\alpha
z^{2}}{4}}\displaystyle\sum_{n=1}^{4} \int_{0}^{S}U_{\alpha}^{(1)}(S,\sigma)P^{\alpha}D_{n}(\sigma)d\sigma\|_{\infty}.$$

Using Proposition ~\ref{ProP}, Equation (~\ref{Compare3}) and the
slow decay of $\beta(\tau)$ we obtain
\begin{equation}
K_{1}\lesssim e^{-\cO S}\|\langle z\rangle^{-3}e^{\frac{\alpha
z^{2}}{4}}\eta(0)\|_{\infty}\lesssim \beta^{2}(T)M_{1}(0).
\end{equation}
By Proposition ~\ref{ProP}, Equations (~\ref{eq:estD1})-
(~\ref{nonlinearity3}) and
$ \int_{0}^{S}e^{-\cO(S-\sigma)}\beta^{2}(\tau(\sigma))d\sigma\lesssim
\beta^{2}(T)$ (see Lemma ~\ref{Bridge}) we have
\begin{equation}\label{M1Ge}
\begin{array}{lll}
K_{2}&\lesssim &\beta^{2}(T)
\{\beta^{\frac{\kappa}{2}}(0)[1+M_{1}(T)A(T)+M_{1}^{2}(T)+M_{1}^{p}(T)]+[M_{2}(T)M_{1}(T)+M_{1}(T)M_{2}^{p-1}(T)]\}.
\end{array}
\end{equation}

 Equation \eqref{NewFun} and the definitions of $S$
and $T$ imply that $\lambda_{1}(t(S))=\lambda(t(T))$, $z=y$,
$\eta(S)=\xi(T)$, and $P^{\alpha}\xi=\xi$, consequently
\begin{equation}\label{eq:ini3}
\|\langle z\rangle^{-3}e^{\frac{\alpha
z^{2}}{4}}P^{\alpha}\eta(S)\|_{\infty}=\|\langle
y\rangle^{-3}e^{\frac{ay^{2}}{4}}\xi(T)\|_{\infty}.
\end{equation}

Collecting the estimates (~\ref{eq:estEsta})-(~\ref{eq:ini3}) and
using the definition of $M_{1}$ in (~\ref{majorants}) we have
$$
\begin{array}{lll}
M_{1}(T)&:=&\displaystyle\sup_{\tau\leq T}\beta^{-2}(\tau)\|\langle
y\rangle^{-3}e^{\frac{ay^{2}}{4}}\xi(\tau)\|_{\infty}\\
&\lesssim &
M_{1}(0)+\beta^{\frac{\kappa}{2}}(0)[1+M_{1}(T)A(T)+M_{1}^{2}(T)+M_{1}^{p}(T)]+M_{2}(T)M_{1}(T)+M_{1}(T)M_{2}^{p-1}(T)
\end{array}
$$
which together with the fact that $T$ is arbitrary implies Equation
(~\ref{M1}).
\begin{flushright}
$\square$
\end{flushright}
\section{Estimate of $M_2$ (Equation \eqref{M2})}\label{SEC:EstM2}
The following lemma is proven similarly to the corresponding parts
of Lemma ~\ref{LM:EstDs} and therefore it is presented without a
proof.
\begin{lemma}\label{LM:keyest}
If $A(\tau), B(\tau)\leq \beta^{-\frac{\kappa}{2}}(\tau)$ and
$b_0\ll 1$ and $D_{n}(\sigma)$, $n=2,3,4$, are the same as in Lemma
~\ref{LM:EstDs}, then
\begin{equation}\label{eq:Remainder0}
\|e^{\frac{\alpha}{4}
z^2}D_{2}(\sigma)\|_{\infty}\lesssim{\beta^{\frac{\kappa}{2}}(\tau
(\sigma))}[\beta^{1/2}(\tau(\sigma))M_{1}(T)+M_{2}(T)];
\end{equation}
\begin{equation}\label{eq:estF0}
\|e^{\frac{\alpha}{4} z^2} D_{3}(\sigma)\|_\infty\lesssim
\beta^{\min\{1,2p-1\}}(\tau(\sigma))[1+ M_1(T)(1+ A(T))+ M_1^2(T)+
M_1^{p}(T)];
\end{equation}
\begin{equation}\label{nonlinearity0}
\|e^{\frac{\alpha }{4}z^{2}}D_{4}(\sigma)\|_{\infty}\lesssim
\beta(\tau(\sigma))M_{1}^{2}(T)+M_{2}^{2}(T)+\beta^{p/2}(\tau(\sigma))M_{1}^{p}(T)+M_{2}^{p}(T).
\end{equation}
\end{lemma}
To estimate $M_2$ it is convenient to treat the $z$-dependent part
of the potential in \eqref{eqn:DefnL0V} as a perturbation.  Let the
operator $L_0$ be the same as in \eqref{eq:eta}. Rewrite
(~\ref{eq:eta}) to have
\begin{equation}
\eta(S)=e^{-(L_{0}+\frac{2p\alpha}{p-1})
S}\eta(0)+ \int_{0}^{S}e^{-(L_{0}+\frac{2p\alpha}{p-1})(S-\sigma)}(V_{2}\eta(\sigma)+\displaystyle\sum_{n=2}^{4}D_{n}(\sigma))d\sigma,
\end{equation} where, recall $S:=\sigma(t(T)),$ $V_{2}$ is the operator given by
$$V_{2}:=\frac{2p\alpha}{p-1+\beta(\tau(\sigma))z^{2}},$$ and the terms $D_{n},\ n=2,3,4,$ are the same as in \eqref{EQ:eta2}. Lemma ~\ref{kernelEst}
implies that
$$\|e^{\frac{\alpha y^{2}}{4}}e^{-(L_{0}+\frac{2p\alpha}{p-1}) s}g\|_{\infty}=e^{-\frac{2p\alpha}{p-1}s}\|e^{\frac{\alpha y^{2}}{4}}e^{-L_{0}s}g\|_{\infty}\lesssim e^{-\frac{2\alpha}{p-1}s}\|e^{\frac{\alpha y^{2}}{4}}g\|_{\infty}$$
for any function $g$ and time $s\geq 0.$ Hence we have
\begin{equation}\label{K123s}
\begin{array}{lll}
\|e^{\frac{\alpha z^{2}}{4}}\eta(S)\|_{\infty}\lesssim
 K_{0}+K_{1}+K_{2}
\end{array}
\end{equation} where the functions $K_{n}$ are given by
$$K_{0}:=e^{-\frac{2\alpha}{p-1} S}\|e^{\frac{\alpha
z^{2}}{4}}\eta(0)\|_{\infty};$$
$$K_{1}:= \int_{0}^{S}e^{-\frac{2\alpha}{p-1}(S-\sigma)}\|e^{\frac{\alpha z^{2}}{4}}V_{2}\eta(\sigma)\|_{\infty}d\sigma,$$
$$K_{2}:=\sum_{n=2}^{4} \int_{0}^{S}e^{-\frac{2\alpha}{p-1}(S-\sigma)}\|e^{\frac{\alpha z^{2}}{4}}D_{n}\|_{\infty}d\sigma.$$

We estimate the $K_{n}$'s, $n=0,1,2.$
\begin{itemize}
\item[(K0)] We start with $K_{0}.$ By (~\ref{Compare0}) and the decay of $e^{-\frac{2\alpha}{p-1} S}$ we
have
\begin{equation}\label{eq:estK0}
K_{0}\lesssim M_{2}(0)+\beta^{1/2}(0)M_{1}(0).
\end{equation}
\item[(K1)] By the definition of
$V_{2}$ we have
$$\|e^{\frac{\alpha z^{2}}{4}}V_{2}\eta(\sigma)\|_{\infty}\lesssim
\|\frac{1}{p-1+\beta(\tau(\sigma))z^{2}}e^{\frac{\alpha
z^{2}}{4}}\eta(\sigma)\|_{\infty}.$$ Moreover by the relation
between $\xi$ and $\eta$ in Equation (~\ref{NewFun}) and Proposition
~\ref{NewTrajectory} we have
$$
\begin{array}{lll}
\displaystyle\max_{0\leq \sigma\leq S} \|e^{\frac{\alpha
z^{2}}{4}}V_{2}\eta(\sigma)\|_{\infty} \lesssim
\displaystyle\max_{T\geq \tau\geq 0}\|\frac{1}{p-1+\beta y^{2}}
e^{\frac{a(\tau)y^{2}}{4}}\xi(\tau)\|_{\infty}.
\end{array}
$$ Using that $D=C/\sqrt{\beta}$ in (~\ref{cutoff}), we find
\begin{equation}\label{EpsiC}
\frac{1}{p-1+\beta y^{2}}\chi_{\geq D}(y)\leq \epsilon(C)
:=\frac{1}{p-1+C^{2}},
\end{equation} which implies
$$\|\frac{1}{p-1+\beta y^{2}}e^{\frac{ay^{2}}{4}}\xi(\tau)\|_{\infty}\leq
\epsilon(C)\|\chi_{\geq
D}e^{\frac{ay^{2}}{4}}\xi(\tau)\|_{\infty}+\|\chi_{\leq
D}e^{\frac{ay^{2}}{4}}\xi(\tau)\|_{\infty}.$$  By the definition of
the function $\chi_{\leq D}$ in Equation (~\ref{cutoff}) we have
that for any $\tau\leq T$, $\chi_{\leq D}\langle
y\rangle^{3}\lesssim \beta^{-3/2}(\tau)$, which implies
$$\|\chi_{\leq D}e^{\frac{ay^{2}}{4}}\xi(\tau)\|_{\infty}\lesssim
\beta^{-3/2}(s)\|\chi_{\leq D}\langle
y\rangle^{-3}e^{\frac{ay^{2}}{4}}\xi(\tau)\|_{\infty}.$$ Collecting
the estimates above, recalling the definitions of $M_{n}, n=1,2,$ in
(~\ref{majorants}), we obtain
\begin{equation}\label{EstK1}
\begin{array}{lll}
K_{1}&\lesssim & \displaystyle\max_{S\geq \sigma\geq
0}\|e^{\frac{\alpha z^{2}}{4}}V_{2}\eta(\sigma)\|_{\infty}
 \int_{0}^{S}e^{-\frac{2\alpha}{p-1}(S-\sigma)}d\sigma \lesssim
\epsilon(C)M_{2}(T)+\beta^{1/2}(0)M_{1}(T).
\end{array}
\end{equation}
\item[(K2)] By the definitions of $D_{n}, \ n=2,3,4,$ and Equations (~\ref{eq:Remainder0})-(~\ref{nonlinearity0}) we have
$$\sum_{n=2}^{4}\|e^{\frac{\alpha z^{2}}{4}}D_{n}(\sigma)\|_\infty\lesssim
\beta^{\frac{\kappa}{2}}(\tau(\sigma))[1+M_{2}(T)+M_{1}(T)A(T)+M_{1}^{2}(T)+M_{1}^{p}(T)]+M_{2}^{2}(T)+M_{2}^{p}(T)$$
and consequently
\begin{equation}\label{EstK2}
K_{2}\lesssim
\beta^{\frac{\kappa}{2}}(0)[1+M_{2}(T)+M_{1}(T)A(T)+M_{1}^{2}(T)+M_{1}^{p}(T)]+M_{2}^{2}(T)+M_{2}^{p}(T).
\end{equation}
\end{itemize}
Collecting the estimates (~\ref{K123s})-(~\ref{EstK2}) we have
\begin{equation}\label{FinalStep}
\begin{array}{lll}
\|e^{\frac{\alpha z^{2}}{4}}\eta(S)\|_{\infty}
&\lesssim& M_{2}(0)+\beta^{1/2}(0)M_{1}(0)+\epsilon(C)M_{2}(T)+\beta^{1/2}(0)M_{1}(T)\\
&
&+\beta^{\frac{\kappa}{2}}(0)[1+M_{2}(T)+M_{1}(T)A(T)+M_{1}^{2}(T)+M_{1}^{p}(T)]+M_{2}^{2}(T)+M_{2}^{p}(T).
\end{array}
\end{equation}
The relation between $\xi$ and $\eta$ in Equation (~\ref{NewFun})
implies
$$
\begin{array}{lll}
\|\chi_{\geq D}e^{\frac{ay^{2}}{4}}\xi(T)\|_{\infty} &\leq&
\|e^{\frac{ay^{2}}{4}}\xi(T)\|_{\infty} =\|e^{\frac{\alpha
z^{2}}{4}}\eta(S)\|_{\infty}
\end{array}
$$
which together with \eqref{FinalStep} gives
$$
\begin{array}{lll}
M_{2}(T)&\lesssim & M_{2}(0)+\beta^{1/2}(0)M_{1}(0)+\epsilon(C)M_{2}(T)+\beta^{1/2}(0)M_{1}(T)+M_{2}^{2}(T)+M_{2}^{p}(T)\\
&
&+\beta^{\frac{\kappa}{2}}(0)[1+M_{2}(T)+M_{1}(T)A(T)+M_{1}^{2}(T)+M_{1}^{p}(T)].
\end{array}
$$
Choosing $C$ so large that $\epsilon(C)$ in Equation (~\ref{EpsiC})
is sufficiently small, we obtain
$$
\begin{array}{lll}
M_{2}(T)&\lesssim& M_{2}(0)+\beta^{1/2}(0)M_{1}(0)+M_{2}^{2}(T)+M_{2}^{p}(T)\\
&
&+\beta^{\frac{\kappa}{2}}(0)[1+M_{2}(T)+M_{1}(T)A(T)+M_{1}^{2}(T)+M_{1}^{p}(T)].
\end{array}
$$
Since $T$ is an arbitrary time, the proof of the estimate \eqref{M2}
for $M_2$ is complete.
\appendix

\section{The Local Well-Posedness of and a Blowup Criterion for \eqref{NLH}}
\label{Appendix:LWP}
In this section we prove the local well-posedness of \eqref{NLH} in
$C([0,T], L^{\infty}).$  The proof is standard and is presented for
the reader's convenience as we did not find it in the literature.
\begin{thm}\label{THM:Local}
Let $u_0\in L^\infty$. For $T=\frac{1}{2}\min[\lb
(2p)^p\|u_0\|_\infty^{p-1}\rb^{-1},1]$ there exists a unique
function $u\in C([0,T],L^\infty)$ satisfying the nonlinear heat
equation \eqref{NLH}.  The solution $u$ depends continuously on the
initial condition $u_0$. Moreover, the solution $u$ satisfies the
estimate
\begin{equation*}
\|u\|_{C([0,T], L^\infty)}\le \max[2^\frac{1}{p}
p\|u_0\|_\infty,2^{\frac{1}{p}}\|u_{0}\|_\infty^{\frac{1}{p}}].
\end{equation*}
Furthermore, either the solution is global in time or blows up in
$L^\infty$ in a finite time.
\end{thm}
\begin{proof}
Using Duhamel's principle, Equation \eqref{NLH} can be written as
the fixed point equation $u=H(u)$, where
\begin{equation}
H(u):=e^{t\p_x^2}u_0+ \int_0^t e^{(t-s)\p_x^2} |u|^{p-1}u(s)\, ds.
\label{NLHFixedPointMap}
\end{equation}
Thus, the proof of existence and uniqueness will be complete if we
can show that the map $H$ has a unique fixed point in the ball
$$B_{R}:=\{u\in X,\ \|u\|_X\leq R\},$$ where $X=:C([0,T],L^\infty)$ and $R:=2\|u_{0}\|_\infty.$
We prove this statement via the contraction mapping principle.

We begin by proving that $H$ is a well-defined map from $B_{R}$ to
$B_{R}$. The estimate
\begin{equation}
\left\|e^{t\p_x^2} u_0\right\|_{X}\le \|u_0\|_\infty
\label{Est:FreeHeat}
\end{equation}
is obtained by using the integral kernel of $e^{t\p_x^2}$,
$e^{t\p_x^2}(x,y)=\frac{1}{\sqrt{\pi t}}e^{-\frac{(x-y)^2}{t}}$,
defined for $t>0$ and its property that $ \int e^{t\p_x^2}(x,y)\,
dy=1$.  Similarly, we find that if $t<T$, then
\begin{equation}
\I{ \int_0^t e^{(t-s)\p_x^2} |u|^{p-1}u(s)\, ds}_X\le T\I{u}_X^p.
\label{Est:NonlinearHeat}
\end{equation}
Estimates \eqref{Est:FreeHeat} and \eqref{Est:NonlinearHeat} imply
that for $T<\infty$, $H:B_R\rightarrow B_R$.

We prove that $H:B_R\rightarrow B_R$ is a strict contraction. Recall
the definition of $T$ in the statement of the theorem. Consider
\begin{equation*}
\|H(u_1)-H(u_2)\|_X\le \I{ \int_0^t \frac{1}{\sqrt{\pi
t}} \int_{-\infty}^\infty e^{-\frac{(x-y)^2}{t}}
|u_1|^{p-1}u_{1}(y,s)-|u_2|^{p-1}u(y,s)|\, dy\, ds}_X.
\end{equation*}
Using that $u_1, u_2\in B_R$, we obtain the estimate $
||u_1|^{p-1}u_{2}-|u_2|^{p-1}u_{2}| \le p |u_1-u_2|R^{p-1}.$ Thus,
\begin{align*}
\|H(u_1)-H(u_2)\|_X&\le p \sup_{[0,T]}\sup_\R
 \int_0^t\frac{1}{\sqrt{\pi t}} \int_{-\infty}^\infty
e^{-\frac{(x-y)^2}{t}}\, dy\, ds \|u_1-u_2\|_X
R^{p-1}\\
&\le p\|u_1-u_2\|_X R^{p-1} T.
\end{align*}
Therefore, if $T<  \frac{1}{2}\min \{ \lb
p^p\|u_0\|_\infty^{p-1}\rb^{-1},1\}$, then $H$ is a strict
contraction in $B_R$. Substituting the choice
$T=\frac{1}{2}\min\{\lb p^p\|u_0\|_\infty^{p-1}\rb^{-1},1\}$ into
the expression for $R$ completes the proof of existence and
uniqueness of $u$ and the estimate on it.

It remains to prove that solution to the initial value problem is
continuous with respect to changes in the initial condition $u_0$.
Let $u$ and $v$ be the solutions with initial conditions $u_0$ and
$v_0$.  We estimate
\begin{align*}
\|u-v\|_X &\leq \|e^{t\p_x^2}(u_0-v_0)\|_{X}+\| \int_0^t
e^{(t-s)\p_x^2}(u^p(s)-v^p(s))\, ds\|_{X}.
\end{align*}
The estimate of these terms proceeds as above (take $u_1=u$ and
$u_2=v$) and if $u,v\in B_R$, then
\begin{equation*}
\|u-v\|_X\leq \|u_0-v_0\|_\infty+\frac{1}{2}\|u-v\|_X.
\end{equation*}
Thus, if $T$ is as above, then $\|u-v\|_X\le 2\|u_0-v_0\|_\infty$
completing the proof of continuity.

Finally, assume $[0,t_*)$ is the maximal interval of existence of
$u$ and $\sup_{0\le t<t_*}\|u(t)\|_\infty:=M<\infty$.  Let
$T:=\frac{1}{2}\min\{ ((2p)^p M^{p-1})^{-1},1 \}$.  Then taking
$u(t_*-\frac{1}{2} T)$ as a new initial condition, we see that the
solution exists in the interval $[0,t_*+\frac{1}{2} T)$, a
contradiction.  This proves the dichotomy claimed in the theorem.
\end{proof}

The theorem below gives a blowup criterion for \eqref{NLH} using the
Lyapunov functional $S(w)$ defined in \eqref{eq:energy}.  Here,
recall, $w(y,s):=(t^*-t)^\frac{1}{p-1}u(x,t)$ with $x=\sqrt{t^*-t}
y$ and $t^*-t=e^{-s}$.  Now we consider $t^*$ as a parameter and
denote $T=t^*-t$.  Then $S(w)=S_T(u)$, where
\begin{equation*}
S_T(u)=T^{\frac{1}{2}\frac{p+3}{p-1}} \int\lb \frac{1}{2}|\nabla
u|^2-\frac{1}{p+1}|u|^{p+1} \rb
\rho(x)\,
dx+\frac{1}{2}\frac{1}{p-1}T^{-\frac{1}{2}\frac{p-5}{p-1}} \int
|u|^2\rho(x )\, dx
\end{equation*}
and $\rho(y)=e^{-\frac{1}{4} y^2}$.

\begin{thm}
Let the initial condition $u_0$ satisfy $S_T(u_0)<0$, modulo a
shift, for some $T>0$. Then \eqref{NLH} blows up in a finite time
$t^*\le T$.
\end{thm}
\begin{proof}
Assume \eqref{NLH} has a solution, $u$, up to time $T$ for an initial
condition $u_0$ as in the theorem.  Let $w$ be as defined in the
paragraph preceding the theorem with $T$ as in the theorem.  The
time derivative of the functional
$I(w):=\frac{1}{2} \int_{-\infty}^\infty w^2(y,s)\rho(y)\, dy$ along
solutions to \eqref{eq:freDir} is
\begin{equation*}
\frac{d}{d s}I(w)=-2 S(w)+\frac{p-1}{p+1} \int_{-\infty}^\infty
|w|^{p+1}\rho\, dy.
\end{equation*}
We use H\"{o}lder's inequality to obtain the estimate
$ \int_{-\infty}^\infty |w|^2\rho\, dy\le
(4\pi)^{\frac{1}{2}\frac{p-1}{p+1}}\lb  \int_{-\infty}^\infty
|w|^{p+1}\rho \rb^\frac{2}{p+1}$.  This and the fact that $S$ is
monotonically decreasing (see \eqref{eqn:9a}) result in the
inequality
\begin{equation*}
\frac{d}{d s}I(w)\ge-2 S(w_0)+\frac{p-1}{p+1}(4\pi)^\frac{1-p}{4}
I(w)^\frac{p+1}{2},
\end{equation*}
and hence if $S(w_0)$ is negative, $I(w)$ blows up in finite time
and therefore so does $w$.  This contradicts our assumption that $u$ exists on $[0,T]$
and, consequently, $w$ exist globally.  To complete the proof, we
write $S(w_0)$ in terms of $S_T(u_0)$.
\end{proof}

\section{Blow-up Dynamics} \label{sec:BlowUpDyn}
In this appendix we investigate the function relation between the
parameters $a$, $b$ and $c$ different from $a=2c-\frac{1}{2}.$

First we observe the following key fact: if $(a,b,c,\xi)$,
$a=f(b,c)$, is a stationary solution to (~\ref{NLH}) satisfying the
estimate $\| \langle y\rangle^{-3} e^{\frac{a}{4} y^2}\xi \|\lesssim
b^2$, then
\begin{equation}
f(0,\frac{1}{2})=\frac{1}{2}. \label{eqn:RedEquilPoint}
\end{equation}
Indeed, if $b=0,$ then the estimate above gives that $\xi=0$ and
therefore $v(y,\tau)=(\frac{2c}{p-1})^{\frac{1}{p-1}}$. Since
$v(y,\tau)$ satisfies (~\ref{eqn:BVNLH}), this implies
(~\ref{eqn:RedEquilPoint}).

In order to simplify our argument, we assume that $f(b,c)$ is of the
form $lc+k$ for some constant $l,k$. By (~\ref{eqn:RedEquilPoint})
we have that $k=\frac{1}{2}-\frac{1}{2}l.$ Thus we have
$a=lc+\frac{1}{2}-\frac{1}{2}l.$

\begin{prop}
For $l>1$,
the different functions $a=lc+\frac{1}{2}-\frac{1}{2}l$ lead to
dynamics equivalent up to rescaling of (~\ref{NLH}).
\end{prop}
\begin{proof}
First we recall that following key points when we prove the case
$a=2c-\frac{1}{2}$, i.e. $l=2$. We decompose the solution of
(~\ref{NLH}) as
\begin{equation}\label{eq:ini}
u_{l=2}(x,t)=\lambda^{\frac{2}{p-1}}(t)[(\frac{2c(\tau)}{p-1+b(\tau)y^{2}})^{\frac{1}{p-1}}+\eta(y,\tau)]
\end{equation}
with $\eta$ satisfying $\|\langle x\rangle^{-3}\eta(x,0)\|=o(b(0))$
and some orthogonality conditions, and $\tau$ and $y$ as defined in
(~\ref{eqn:split2}). And for any $l$ we define
\begin{equation}\label{defna}
a(t(\tau)):=\lambda^{-3}(t)\frac{d}{dt}\lambda(t).
\end{equation} We
require $2c(0)=c_{l=2}(0)=1-\frac{2}{p-1}b(0)+O(b^{2}(0))$. Using
Equations (~\ref{eqn:bParameter}) and (~\ref{eqn:cParameter}) we get
that $2c_{l=2}(\tau)=1-\frac{2}{p-1}b(\tau)+O(b^{2})$ and
$b(\tau)\rightarrow 0^+$, $\frac{dc_{l=2}(\tau)}{d\tau}=O(b^{3})$.
On the other hand we have that if
$2c(0)=2c_{l}(0)=1+\frac{2}{(1-l)(p-1)}b(0)+O(b^{2}(0))$ and
$\|\langle x\rangle^{-3}\eta(x,0)\|=o(b(0))$, we fix the function as
\begin{equation}\label{fixfun}
a=lc_{l}+\frac{1}{2}-\frac{1}{2}l,
\end{equation} after going
through the same procedure we prove that
$2c_{l}(\tau)=1+\frac{2}{(1-l)(p-1)}b(\tau)+O(b^{2})$,
$b(\tau)\rightarrow 0^+$, $\frac{d}{d\tau}c_{l}(\tau)=O(b^{3}).$ The
two equations are related to each other in the following sense.

If $c(0)$ in (~\ref{eq:ini}) satisfies the condition that
$c(0)=c_{l}(0)=1+\frac{2}{(1-l)(p-1)}b(0)+O(b^{2})$ for $l>1$ then
we rewrite
$$u_{l=2}(y,\tau)=\lambda_{1}^{\frac{2}{p-1}}(t)[(\frac{2c_{1}(\tau)}{p-1+\beta(\tau)y_{1}^{2}})^{\frac{1}{p-1}}+\eta_{2}(y_{1},\tau)]$$ with
$\lambda_{1}(t):=\lambda(t)\sqrt{\frac{c_{l=2}(\tau(t))}{c_{l_{0}}(\tau(t))}}$,
$y_{1}:=\lambda_{1}(t)x$ and
$\beta(\tau):=b(\tau)\frac{c_{l_{0}}(\tau)}{c_{l=2}(\tau)}$ and
$\eta_{2}$ from $\eta(y,\tau)=o(b).$ We compute to get
$$a_{1}:=\lambda_{1}^{-3}(t)\frac{d}{dt}\lambda_{1}(t)=a_{l=2}(\frac{c_{l=2}(\tau(t))}{c_{l}(\tau(t))})^{2}+O(b^{3})=a_{l}+O(b^{2})$$
$$\frac{d}{d\tau}\beta=-\frac{4p}{(p-1)^{2}}\beta^{2}+O(b^{3})$$
thus $a_{1}=lc_{l}+\frac{1}{2}-\frac{1}{2}l+O(b^{2})$ which is
consistent with (~\ref{defna}) and (~\ref{fixfun}) (the remainder
$O(b^{2})$ in the function of $a_{1}$ can be erased by adding some
correction on $c_{l}$). Thus the case $l=2$ can be transformed into
the other $l>1$ cases. By similar argument we prove that all these
are equivalent.
\end{proof}
Now we remark on the dynamics of the parameters $a$, $b$ and $c$
described by Equations \eqref{eqn:bParameter} and
\eqref{eqn:cParameter} if we neglect the remainder terms determined
by the fluctuations $\xi$. In other words we consider the truncated
dynamical system for the parameters $b$ and $c$ which reads
\begin{align}
b_\tau&=-\frac{2}{p-1}\lb 1+\frac{2 p}{p-1}\rb b^2+2(c-a)b+O(b^{3})
,\label{eqn:bParameterXiZero}\\
c_\tau&=2 c(c-a)-\frac{2}{p-1}b
c+O(b^{3}).\label{eqn:cParameterXiZero}
\end{align}
A simple computation shows that if $a=lc+\frac{1}{2}-\frac{1}{2}l$
and $l>1$, then the point $(b,c)=(0,\frac{1}{2})$ is marginally
stable for (~\ref{eqn:bParameterXiZero}) and
(~\ref{eqn:cParameterXiZero}).

\section{Spectrum of the Linear Operator ${\mathcal L}_{a b c}$}
\label{Appendix:SpectrumLinear}
We assume that the $|a_\tau|$ term is negligible in comparison with
$a$ and consider the operator $\tilde{{\cal L}}_{a b c}$, which
differs from ${\cal L}_{a b c}$ by the term $\frac{1}{4}a_\tau y^2$:
\begin{equation*}
\widetilde{\cal L}_{a b c}:=-\p_y^2+\frac{1}{4} a^2
y^2-\frac{a}{2}+\frac{2 a}{p-1}-\frac{2 p c}{p-1+ b y^2}.
\end{equation*}
Due to the quadratic term $\frac{1}{4} a y^2$, the operator
$\widetilde{\L}_{a b c}$ has a purely discrete spectrum. We can
obtain a better understanding of its eigenvalues by comparing it to
the harmonic oscillator
\begin{equation}
\L_0:=-\p_y^2+\frac{1}{4} a^2 y^2-\frac{a}{2}.
\end{equation}
Then $\L_{0}+\frac{2}{p-1}(a-pc)$ and $\L_{0}+\frac{2a}{p-1}$
approximate $\tilde{\cal L}_{a b c}$ near zero and at infinity,
respectively.  The spectrum of the operator $\L_{0}$ is
\begin{equation}
\sigma\lb{\mathcal L}_0\rb=\left\{n a |\ n=0,1,2,\ldots\right\}.
\end{equation}
The first three normalized eigenvectors of ${\mathcal L}_0$, which
are used in the main part of the paper, are
\begin{equation}\label{eq:eigenvectors}
\phi_{0a}:=\lb\frac{a}{2\pi}\rb^\frac{1}{4} e^{-\frac{a}{4}y^2},\
\phi_{1a}:=\lb\frac{a}{2\pi}\rb^{\frac{1}{4}}\sqrt{a}y
e^{-\frac{a}{4}y^2},
\phi_{2a}:=\lb\frac{a}{8\pi}\rb^{\frac{1}{4}}(1-a
y^2)e^{-\frac{a}{4}y^2}.
\end{equation}

\begin{prop}
If $p>1$, $c\ge0$ and $b\ge0$, then the eigenvalues $\lambda_n$ of
$\widetilde{\mathcal L}_{a b c}$ satisfy the bounds
\begin{equation}
n a+\frac{2 a}{p-1}\ge\lambda_n\ge n a+\frac{2}{p-1}(a-p c).
\label{eqn:EigenvalueBounds}
\end{equation}
\end{prop}
\begin{proof}
First we show that
\begin{equation}
\L_0+\frac{2a}{p-1}>\widetilde{\L}_{a b c}>
\L_0+\frac{2}{p-1}(a-pc). \label{eqn:74}
\end{equation}
Since $p>1$, $b\ge 0$ and $c\ge 0$, $0<\frac{2 p c}{p-1+b
y^2}\le\frac{2 p c}{p-1}$, and hence \eqref{eqn:74}.  The $n$-th
eigenvalue of $\widetilde{\L}_{a b c}$ (starting from $n=0$) is by
the MinMax principle
\begin{equation}
\lambda_n=\sup_{\dim X=n}\inf_{\{\psi\in
X^\bot|\|\psi\|=1\}}\ip{\psi}{\widetilde{\L}_{a b c}\psi}.
\end{equation}
Using the inequality $\langle \psi,\widetilde{\L}_{a b
c}\psi\rangle\ge\ip{\psi}{\L_0\psi}+\frac{2}{p-1}(a-pc)\langle \psi,\psi\rangle$ and the characterization of the
spectrum of $\L_0$ we obtain
\begin{equation}
\lambda_n\ge\sup_{\dim X=n}\inf_{\{\psi\in
X^\bot|\|\psi\|=1\}}\ip{\psi}{\L_0\psi}+\frac{2}{p-1}(a-pc)=n a+\frac{2}{p-1}(a-pc)
\end{equation}
and similarly for the upper bound.
\end{proof}

Equation \eqref{eqn:cParameter} and the relation $a=2c-\frac{1}{2}$
suggests that $c=a+\O{b}$ where $b$ is small.  In this case Equation
\eqref{eqn:EigenvalueBounds} shows that the operator
$\widetilde{\cal L}_{a b c}$ has at most three non-positive
eigenvalues. The second eigenvalue corresponds to an odd
eigenfunction and therefore drops out if we assume that the initial
condition $u_0(x)$ is even (so that $x_0=0$, otherwise one has to
use the parameter $x_0$). The two parameters $b$ and $c$ are chosen
so that the fluctuation $\xi$ is orthogonal to the other two
eigenfunctions. Hence on the space of $\xi$'s the linear operator
$\widetilde{\cal L}_{a b c}$ has strictly positive spectrum.

\section{Proof of the Feynmann-Kac Formula}\label{Sec:Trotter}
In this appendix we present, for the reader's convenience, a proof
of the Feynman-Kac formula (~\ref{eq:FK1})-(~\ref{FK2}) and the
estimate \eqref{eqn:93a} (cf. \cite{BrKu}). For stochastic calculus
proofs of similar formulae see \cite{Du, GlJa, Hida, KaSh, Simon}.

Let
$L_{0}:=-\partial_{y}^{2}+\frac{\alpha^{2}}{4}y^{2}-\frac{\alpha}{2}$
and $L:=L_{0}+V$ where $V$ is a multiplication operator by a
function $V(y,\tau)$, which is bounded and Lipschitz continuous in
$\tau$. Let $U(\tau,\sigma)$ and $U_{0}(\tau,\sigma)$ be the
propagators generated by the operators $-L$ and $-L_{0},$
respectively. The integral kernels of these operators will be
denoted by $U(\tau,\sigma)(x,y)$ and $U_0(\tau,\sigma)(x,y)$.
\begin{thm}\label{THM:trotter}
The integral kernel of $U(\tau,\sigma)$ can be represented as
\begin{equation}
U(\tau,\sigma)(x,y)=U_{0}(\tau,\sigma)(x,y) \int
e^{  \int_{\sigma}^{\tau}V(\omega_{0}(s)+\omega(s),s)ds}d\mu(\omega)
\label{eqn:BNeg1}\end{equation} where $d\mu(\omega)$ is a
probability measure (more precisely, a conditional harmonic
oscillator, or Ornstein-Uhlenbeck, probability measure) on the
continuous paths $\omega: [\sigma,\tau]\rightarrow\R$ with
$\omega(\sigma)=\omega(\tau)=0$, and $\omega_{0}(\cdot)$ is the path
defined as
\begin{equation}
\omega_{0}(s)=e^{\alpha (\tau-s)}\frac{e^{2\alpha\sigma}-e^{2\alpha
s}}{e^{2\alpha\sigma}-e^{2\alpha \tau}}x+e^{\alpha (\sigma-s)}\frac{e^{2\alpha\tau}-e^{2\alpha s}}{e^{2\alpha \tau}-e^{2\alpha\sigma}}y.
\label{eqn:B0}
\end{equation}
\end{thm}
\begin{remark}
\label{remark:OrnsteinUhlenbeck} $d\mu(\omega)$ is the Gaussian
measure with mean zero and covariance
$(-\p_s^2+\alpha^2)^{-1}$, normalized to 1. The path
$\omega_0(s)$ solves the boundary value problem
\begin{equation}
(-\p_s^2+\alpha^2)\omega_0=0\ \mbox{with}\
\omega(\sigma)=y\ \mbox{and}\ \omega(\tau)=x. \label{eqn:133a}
\end{equation}
Below we will also deal with the normalized Gaussian measure
$d\mu_{x y}(\omega)$ with mean $\omega_0(s)$ and covariance
$(-\p_s^2+\alpha^2)^{-1}$. This is a conditional
Ornstein-Uhlenbeck probability measure on continuous paths
$\omega:[\sigma,\tau]\rightarrow\R$ with $\omega(\sigma)=y$ and
$\omega(\tau)=x$ (see e.g. \cite{GlJa, Hida, Simon}).
\end{remark} Now, assume in addition that the function $V(y,\tau)$ satisfies the
estimates \begin{equation}\label{eq:diffeV} V\leq 0\ \text{and}\
|\partial_{y}V(y,\tau)|\lesssim
\beta^{-\frac{1}{2}}(\tau)\end{equation} where $\beta(\tau)$ is a
positive function. Then Theorem ~\ref{THM:trotter} implies Equation
(~\ref{eqn:93a}) by the following corollary.
\begin{cor}\label{cor:diff} Under (~\ref{eq:diffeV}),
$$|\partial_{y}  \int e^{  \int_{\sigma}^{\tau}V(\omega_{0}(s)+\omega(s),s)ds}d\mu(\omega)|\lesssim |\tau-\sigma|\sup _{\sigma\leq s\leq \tau} \beta^{1/2}(\tau) $$
\end{cor}
\begin{proof}
By Fubini's theorem $$\partial_{y}  \int
e^{  \int_{\sigma}^{\tau}V(\omega_{0}(s)+\omega(s),s)ds}d\mu(\omega)= \int
\partial_{y}[ \int_{0}^{\tau}V(\omega_{0}(s)+\omega(s),s)ds]e^{ \int_{\sigma}^{\tau}V(\omega_{0}(s)+\omega(s),s)ds}d\mu(\omega)$$
Equation (~\ref{eq:diffeV}) implies$$
|\partial_{y} \int_{\sigma}^{\tau}V(\omega_{0}(s)+\omega(s),s)ds|\leq
|\tau-\sigma| \sup _{\sigma\leq s\leq \tau} \beta^{1/2}(\tau) |,\
\text{and}\
e^{ \int_{\sigma}^{\tau}V(\omega_{0}(s)+\omega(s),s)ds}\leq 1.$$ Thus
$$|\partial_{y}  \int
e^{ \int_{\sigma}^{\tau}V(\omega_{0}(s)+\omega(s),s)ds}d\mu(\omega)|\lesssim
|\tau-\sigma| \sup _{\sigma\leq s\leq \tau} \beta^{1/2}(\tau) | \int
d\mu(\omega)=|\tau-\sigma| \sup _{\sigma\leq s\leq \tau}
\beta^{1/2}(\tau) |
$$ to complete the proof.
\end{proof}

\begin{proof}[Proof of Theorem \ref{THM:trotter}]
We begin with the following extension of the Ornstein-Uhlenbeck
process-based Feynman-Kac formula to time-dependent potentials:
\begin{equation}
U(\tau,\sigma)(x,y)=U_0(\tau,\sigma)(x,y) \int e^{- \int_\sigma^\tau
V(\omega(s),s)\, ds}d\mu_{x y}(\omega). \label{eqn:B2}
\end{equation}
where $d\mu_{x y}(w)$ is the conditional Ornstein-Uhlenbeck
probability measure described in Remark
\ref{remark:OrnsteinUhlenbeck} above. This formula can be proven in
the same way as the one for time independent potentials (see
\cite{GlJa}, Equation (3.2.8)), i.e. by using the Kato-Trotter
formula and evaluation of Gaussian measures on cylindrical sets.
Since its proof contains a slight technical wrinkle, for the
reader's convenience we present it below.

Now changing the variable of integration in \eqref{eqn:B2} as
$\omega=\omega_0+\tilde{\omega}$, where $\tilde{\omega}(s)$ is a
continuous path with boundary conditions
$\tilde{\omega}(\sigma)=\tilde{\omega}(\tau)=0$, using the
translational change of variables formula $ \int f(\omega)\,
d\mu_{xy}(\omega)= \int f(\omega_0+\tilde{\omega})\,
d\mu(\tilde{\omega})$, which can be proven by taking
$f(\omega)=e^{i\langle\omega,\zeta\rangle}$ and using
\eqref{eqn:133a} (see \cite{GlJa}, Equation (9.1.27)) and omitting
the tilde over $\omega$ we arrive at \eqref{eqn:BNeg1}.
\end{proof}
There are at least three standard ways to prove \eqref{eqn:B2}: by
using the Kato-Trotter formula, by expanding both sides of the
equation in $V$ and comparing the resulting series term by term and
by using Ito's calculus (see \cite{KaSh, Simon,RSII,GlJa}). The
first two proofs are elementary but involve tedious estimates while
the third proof is based on a fair amount of stochastic calculus.
For the reader's convenience, we present the first elementary proof
of \eqref{eqn:B2}.

Before starting proving \eqref{eqn:B2} we establish an auxiliary
result. We define the operator $\mathcal{K}$ as
\begin{equation}\label{eqn:defK}
\mathcal{K}(\sigma,\delta):= \int_{0}^{\delta}U_{0}(\sigma+\delta,\sigma+s)V(\sigma+s,\cdot)U_{0}(\sigma+s,\sigma)ds-U_{0}(\sigma+\delta,\sigma) \int_{0}^{\delta}
V(\sigma+s,\cdot)ds
\end{equation}

\begin{lemma} For any $\sigma\in [0,\tau]$ and $\xi\in\mathcal{C}_{0}^{\infty}$ we
have, as $\delta\rightarrow 0^+$,
\begin{equation}\label{eqn:appro}
\sup_{0\leq \sigma\leq
\tau}\|\frac{1}{\delta}\mathcal{K}(\sigma,\delta)U(\sigma,0)\xi\|_{2}\rightarrow
0.
\end{equation}
\end{lemma}
\begin{proof}
If the potential term, $V$, is independent of $\tau$, then the proof
is standard (see, e.g. \cite{RSII}). We use the property that the
function $V$ is Lipschitz continuous in time $\tau$ to prove
(~\ref{eqn:appro}). The operator $\mathcal{K}$ can be further
decomposed as
$$\mathcal{K}(\sigma,\delta)=\mathcal{K}_{1}(\sigma,\delta)+\mathcal{K}_{2}(\sigma,\delta)$$ with
$$\mathcal{K}_{1}(\sigma,\delta):= \int_{0}^{\delta}U_{0}(\sigma+\delta,\sigma+s)
V(\sigma,\cdot)U_{0}(\sigma+s,\sigma)ds-\delta
U_{0}(\sigma+\delta,\sigma) V(\sigma,\cdot)$$ and
$$\mathcal{K}_{2}(\sigma,
\delta):= \int_{0}^{\delta}U_{0}(\sigma+\delta,\sigma+s)
[V(\sigma+s,\cdot)-V(\sigma,\cdot)]U_{0}(\sigma+s,\sigma)ds-U_{0}(\sigma+\delta,\sigma)
 \int_{0}^{\delta} [V(\sigma+s,\cdot)-V(\sigma,\cdot)]ds.$$

Since $U_{0}(\tau,\sigma)$ are uniformly $L^{2}$-bounded and $V$ is
bounded, we have $U(\tau,\sigma)$ is uniformly $L^2$-bounded. This
together with the fact that the function $V(\tau,y)$ is Lipschitz
continuous in $\tau$ implies that
$$\|\mathcal{K}_{2}(\sigma,\delta)\|_{L^{2}\rightarrow L^{2}}\lesssim 2 \int_{0}^{\delta}sds=\delta^{2}.$$

We rewrite $\mathcal{K}_{1}(\sigma,\delta)$ as
$$\mathcal{K}_{1}(\sigma,\delta)= \int_{0}^{\delta} U_{0}
(\sigma+\delta,\sigma+s)\{V(\sigma,\cdot)[U_{0}(\sigma+s,\sigma)-1]-
[U_{0}(\sigma+s,\sigma)-1]V(\sigma,\cdot)\}ds.$$  Let
$\xi(\sigma)=U(\sigma,0)\xi$.  We claim that for a fixed $\sigma\in
[0,\tau]$,
\begin{equation}
\|\mathcal{K}_{1}(\sigma,\delta)\xi(\sigma)\|_{2}=o(\delta).
\label{eqn:142a}
\end{equation}
 Indeed, the fact
$\xi_{0}\in \mathcal{C}_{0}^{\infty}$ implies that
$L_{0}\xi(\sigma),\ L_{0}V(\sigma)\xi(\sigma)\in L^{2}.$
Consequently (see \cite{RSI})
$$\lim_{s\rightarrow
0^+}\frac{(U_{0}(\sigma+s,\sigma)-1)g}{s}\rightarrow L_{0}g, $$ for
$g=\xi(\sigma)\ \text{or}\ V(\sigma,y)\xi(\sigma)$ which implies our
claim. Since the set of functions $\{\xi(\sigma)|\sigma\in
[0,\tau]\}\subset L_{0}L^{2}$ is compact and
$\|\frac{1}{\delta}K_{1}(\sigma,\delta)\|_{L^2\rightarrow L^2}$ is
uniformly bounded, we have \eqref{eqn:142a} as $\delta\rightarrow 0$
uniformly in $\sigma\in [0,\tau]$.

Collecting the estimates on the operators $\mathcal{K}_{i},\ i=1,2$,
we arrive at (~\ref{eqn:appro}).
\end{proof}
\begin{lemma}
Equation \eqref{eqn:B2} holds.
\end{lemma}
\begin{proof} In order to simplify our notation, in the proof that follows we assume,
without losing generality, that $\sigma=0$.  We divide the proof
into two parts. First we prove that for any fixed $\xi\in
\mathcal{C}_{0}^{\infty}$ the following Kato-Trotter type formula holds
\begin{equation}\label{eqn:trotter}
U(\tau,0)\xi=\lim_{n\rightarrow \infty}\prod_{0\leq k\leq
n-1}U_{0}(\frac{k+1}{n}\tau,
\frac{k}{n}\tau)e^{ \int_{\frac{k\tau}{n}}^{\frac{(k+1)\tau}{n}}V(y,s)ds}\xi
\end{equation} in the $L^{2}$ space.
We start with the formula
$$
\begin{array}{lll}
& &U(\tau,0)-\displaystyle\prod_{0\leq k\leq n-1}U_{0}(\frac{k+1}{n}\tau, \frac{k}{n}\tau)e^{ \int_{\frac{k\tau}{n}}^{\frac{(k+1)\tau}{n}}V(y,s)ds}\\
&=&\displaystyle\prod_{0\leq k\leq n-1}U(\frac{k+1}{n}\tau, \frac{k}{n}\tau)-\prod_{0\leq k\leq n-1}U_{0}(\frac{k+1}{n}\tau, \frac{k}{n}\tau)e^{ \int_{\frac{k\tau}{n}}^{\frac{(k+1)\tau}{n}}V(y,s)ds}\\
&=&\displaystyle\sum_{0\leq j\leq n}\prod_{j\leq k\leq
n-1}U_{0}(\frac{k+1}{n}\tau,
\frac{k}{n}\tau)e^{ \int_{\frac{k\tau}{n}}^{\frac{(k+1)\tau}{n}}V(y,s)ds}A_{j}U(\frac{j}{n}\tau,
0)
\end{array}
$$ with the operator $$A_{j}:=U_{0}(\frac{j+1}{n}\tau, \frac{j}{n}\tau)e^{ \int_{\frac{j\tau}{n}}^{\frac{(j+1)\tau}{n}}V(y,s)ds}-U(\frac{j+1}{n}\tau, \frac{j}{n}\tau).$$

We observe that $\|U_{0}(\tau,\sigma)\|_{L^{2}\rightarrow L^{2}}\leq 1$, and moreover by the boundness of $V,$ the operator $U(\tau,\sigma)$
is uniformly bounded in $\tau$ and $\sigma$ in any compact set. Consequently
\begin{equation}\label{eq:b4}
\begin{array}{lll}
& &\|[U(\tau,0)-\displaystyle\prod_{0\leq k\leq n-1}U_{0}(\frac{k+1}{n}\tau, \frac{k}{n}\tau)e^{ \int_{\frac{k\tau}{n}}^{\frac{(k+1)\tau}{n}}V(y,s)ds}]\xi\|_{2}\\
&\leq &\displaystyle\max_{j}n\|\displaystyle\prod_{j\leq k\leq n-1}U_{0}(\frac{k+1}{n}\tau, \frac{k}{n}\tau)e^{ \int_{\frac{k\tau}{n}}^{\frac{(k+1)\tau}{n}}V(y,s)ds}A_{j}U(\frac{j}{n}\tau, 0)\xi\|_{2}\\
&\lesssim
&n\displaystyle\max_{j}\|A_{j}+\mathcal{K}(\frac{k}{n}\tau,
\frac{1}{n}\tau)\|_{L^{2}\rightarrow
L^{2}}+\displaystyle\max_{j}n\|\mathcal{K}(\frac{j}{n}\tau,\frac{1}{n}\tau)U(\frac{j}{n},0)\xi\|_{2}
\end{array}
\end{equation}
where, recall the definition of $\mathcal{K}$ from
(~\ref{eqn:defK}). Now we claim that
\begin{equation}\label{eqn:assum}
\|A_{j}+\mathcal{K}(\frac{k}{n}\tau,
\frac{1}{n}\tau)\|_{L^2\rightarrow L^2}\lesssim \frac{1}{n^{2}}.
\end{equation}

Indeed, by Duhamel's principle we have
$$U(\frac{j+1}{n}\tau,\frac{j}{n}\tau)=U_{0}(\frac{j+1}{n}\tau,\frac{j}{n}\tau)+ \int_{0}^{\frac{1}{n}\tau}U_{0}(\frac{j+1}{n}\tau,s)V(y,s)U(s,\frac{j}{n}\tau)ds.$$
Iterating this equation on $U(s,\frac{k}{n}\tau)$ and using the fact that $U(s,t)$ is uniformly bounded if $s,t$ is on a compact set, we obtain
$$\|U(\frac{j+1}{n}\tau,\frac{j}{n}\tau)-U_{0}(\frac{j+1}{n}\tau,\frac{j}{n}\tau)- \int_{0}^{\frac{1}{n}\tau}U_{0}(\frac{j+1}{n}\tau,s)V(y,s)U_{0}(s,\frac{j}{n}\tau)ds\|_{L^{2}\rightarrow L^{2}}\lesssim \frac{1}{n^{2}}.$$
On the other hand we expand
$e^{ \int_{\frac{j\tau}{n}}^{\frac{(j+1)\tau}{n}}V(y,s)ds}$ and use the fact that $V$ is bounded to get
$$\|U_{0}(\frac{j+1}{n}\tau,\frac{j}{n}\tau)e^{ \int_{\frac{j\tau}{n}}^{\frac{(j+1)\tau}{n}}V(y,s)ds}-U_{0}(\frac{j+1}{n}\tau,\frac{j}{n}\tau)-U_{0}(\frac{j+1}{n}\tau,\frac{j}{n}\tau) \int_{\frac{j\tau}{n}}^{\frac{(j+1)\tau}{n}}V(y,s)ds\|_{L^{2}\rightarrow L^{2}}\lesssim \frac{1}{n^{2}}.$$
By the definition of $\mathcal{K}$ and $A_{j}$ we complete the proof
of (~\ref{eqn:assum}). Equations (~\ref{eqn:appro}), (~\ref{eq:b4})
and (~\ref{eqn:assum}) imply (~\ref{eqn:trotter}). This completes
the first step.

In the second step we compute the integral kernel, $G_{n}(x,y)$, of
the operator $$G_{n}:=\displaystyle\prod_{0\leq k\leq
n-1}U_{0}(\frac{k+1}{n}\tau,
\frac{k}{n}\tau)e^{ \int_{\frac{k\tau}{n}}^{\frac{(k+1)\tau}{n}}V(\cdot,s)ds}$$
in (~\ref{eqn:trotter}). By the definition, $G_{n}(x,y)$ can be
written as
\begin{equation}
G_{n}(x,y)= \int\cdot\cdot\cdot \int \prod_{0\leq k\leq
n-1}U_{\frac{\tau}{n}}(x_{k+1},x_{k})e^{ \int_{\frac{k\tau}{n}}^{\frac{(k+1)\tau}{n}}V(x_{k},s)ds}dx_{1}\cdot\cdot\cdot
dx_{n-1} \label{eqn:TrotterAsterik}
\end{equation}
 with $x_{n}:=x,\ x_{0}:=y$ and $U_\tau(x,y)\equiv U_0(0,\tau)(x,y)$
 is the integral kernel of the operator $U_0(\tau,0)=e^{-L_0\tau}$.
We rewrite \eqref{eqn:TrotterAsterik} as
\begin{equation}
 G_n(x,y)=U_\tau(x,y) \int e^{\sum_{k=0}^{n-1} \int_{\frac{k\tau}{n}}
 ^{\frac{(k+1)\tau}{n}} V(x_k,s)\, ds}\, d\mu_n(x_1,\ldots, x_n),
\label{eqn:TrotterDoubAsterik1}
\end{equation}
where
\begin{equation*}
d\mu_n(x_1,\ldots, x_n):=\frac{\prod_{0\le k\le
n-1}U_{\frac{\tau}{n}}(x_{k+1},x_k)}{U_\tau(x,y)}dx_1\ldots
dx_{k-1}.
\end{equation*}
Since $G_n(x,y)|_{V=0}=U_\tau(x,y)$ we have that $ \int
d\mu_n(x_1,\ldots,x_n)=1$.  Let
$\Delta:=\Delta_1\times\ldots\times\Delta_n$, where $\Delta_j$ is an
interval in $\R$.  Define a cylinderical set
\begin{equation*}
P^n_\Delta:=\{ \omega:[0,\tau]\rightarrow\R\ |\ \omega(0)=y,\
\omega(\tau)=x,\ \omega(k\tau/n)\in \Delta_k,\ 1\le k\le n-1 \}.
\end{equation*}
By the definition of the measure $d\mu_{xy}(\omega)$, we have
$\mu_{xy}(P^n_\Delta)= \int_{\Delta} d\mu_n(x_1,\ldots, x_n)$.  Thus,
we can rewrite \eqref{eqn:TrotterDoubAsterik1} as
\begin{equation}
 G_n(x,y)=U_\tau(x,y) \int e^{\sum_{k=0}^{n-1}
  \int_{\frac{k\tau}{n}}^{\frac{(k+1)\tau}{n}} V(\omega(\frac{k\tau}{n}),s)\, ds}\, d\mu_{xy}(\omega),
\label{eqn:TrotterDoubAsterik2}
\end{equation}
By the dominated convergence theorem the integral on the right hand
side of \eqref{eqn:TrotterDoubAsterik2} converges in the sense of
distributions as $n\rightarrow\infty$ to the integral on the right
hand side of \eqref{eqn:B2}.  Since the left hand side of
\eqref{eqn:TrotterDoubAsterik2} converges to the left hand side of
\eqref{eqn:B2}, also in the sense of distributions (which follows
from the fact that $G_n$ converges in the operator norm on $L^2$ to
$U(\tau,\sigma)$), \eqref{eqn:B2} follows.
\end{proof}

Note that on the level of finite dimensional approximations the
change of variables formula can be derived as follows.  It is
tedious, but not hard, to prove that
$$\prod_{0\leq k\leq n-1}U_{n}(x_{k+1},x_{k})=e^{-\alpha\frac{(x-e^{-\alpha \tau}y)^{2}}{2(1-e^{-2\alpha \tau})}}\prod_{0\leq k\leq n-1}U_{n}(y_{k+1},y_{k})$$
with $y_{k}:=x_{k}-\omega_{0}(\frac{k}{n}\tau)$. By the definition
of $\omega_{0}(s)$ and the relations $x_{0}=y$ and $x_{n}=x$ we have
\begin{equation}\label{eq:gnzy}
G_{n}(x,y)=U_\tau(x,y)G^{(1)}_{n}(x,y)
\end{equation}
where
\begin{equation}\label{eq:g1xy}
G^{(1)}_{n}(x,y):=\frac{1}{4\pi \sqrt{\alpha}(1-e^{-2\alpha
\tau})} \int\cdot\cdot\cdot \int\prod_{0\leq k\leq
n-1}U_{n}(y_{k+1},y_{k})e^{ \int_{\frac{k\tau}{n}}^{\frac{(k+1)\tau}{n}}V(y_{k}+\omega_{0}(\frac{k\tau}{n}),s)ds}dy_{1}\cdot\cdot\cdot
dy_{k-1}.\end{equation} Since $\displaystyle\lim_{n\rightarrow
\infty}G_{n}\xi$ exists by (~\ref{eqn:appro}), we have
$\displaystyle\lim_{n\rightarrow\infty}G^{(1)}_{n}\xi$ (in the weak
limit) exists also. As shown in \cite{GlJa},
$\displaystyle\lim_{n\rightarrow\infty}G^{(1)}_{n}= \int
e^{ \int_{0}^{\tau}V(\omega_{0}(s)+\omega(s),s)ds}d\mu(\omega)$ with
$d\mu$ being the (conditional) Ornstein-Uhlenbeck measure on the set
of path from $0$ to $0.$ This completes the derivation of the change
of variables formula.
\begin{remark}
In fact, Equations (~\ref{eqn:trotter}), (~\ref{eq:gnzy}) and
(~\ref{eq:g1xy}) suffice to prove the estimate in Corollary
~\ref{cor:diff}.
\end{remark}



\end{document}